\documentclass{amsart}
\usepackage{amssymb,stmaryrd}
\usepackage{amsfonts}
\usepackage{amstext}
\usepackage{algorithmic}
\usepackage{algorithm}
\usepackage{graphicx}
\usepackage{epstopdf}
\usepackage[all]{xy}

\numberwithin{equation}{section}
\newtheorem{theorem}{Theorem}[section]
\newtheorem{lemma}[theorem]{Lemma}
\newtheorem{proposition}[theorem]{Proposition}
\newtheorem{corollary}[theorem]{Corollary}

\theoremstyle{definition}

\theoremstyle{remark}

\newcommand{\R}{{\mathbb{R}}}

\newcommand{\Q}{{\mathbb{Q}}}
\newcommand{\F}{{\mathbb{F}}}

\newcommand{\<}{{\langle}}
\renewcommand{\(}{{(}}
\renewcommand{\)}{{)}}
\renewcommand{\>}{{\rangle}}

\newcommand{\CE}{{\mathcal{E}}}

\newcommand{\CL}{{\mathcal{L}}}

\newcommand{\tens}{\otimes}
\newcommand{\btens}{{\bar\otimes}}
\newcommand{\htens}{{\hat\otimes}}
\newcommand{\bnabla}{{\bar\nabla}}
\newcommand{\bd}{{\bar\extd}}
\newcommand{\bDelta}{{\bar\Delta}}

\newcommand{\id}{\rm id}

\newcommand{\extd}{{\rm d}}
\newcommand{\del}{{\partial}}
\newcommand{\eps}{\epsilon}

\begin{document}

\title{Almost commutative Riemannian geometry:  wave operators}
\keywords{noncommutative geometry, quantum groups, quantum gravity, Riemannian manifold, Ricci curvature, Schwarzschild black hole, gravitational redshift}

\subjclass[2000]{Primary 81R50, 58B32, 83C57}

\author{Shahn Majid}
\address{Queen Mary, University of London\\
School of Mathematics, Mile End Rd, London E1 4NS, UK}

\email{s.majid@qmul.ac.uk}

\thanks{The author was supported by a Leverhulme Senior Research Fellowship}

\begin{abstract}
Associated to any (pseudo)-Riemannian manifold $M$ of dimension $n$ is an $n+1$-dimensional noncommutative differential structure $(\Omega^1,\extd)$ on the manifold, with the extra dimension encoding the classical Laplacian as a noncommutative `vector field'. We use the classical connection, Ricci tensor and Hodge Laplacian to construct $(\Omega^2,\extd)$ and a natural noncommutative torsion free connection $(\nabla,\sigma)$ on $\Omega^1$. We show that its generalised braiding $\sigma:\Omega^1\tens\Omega^1\to \Omega^1\tens\Omega^1$ obeys the quantum Yang-Baxter or braid relations only when the original $M$ is flat, i.e their failure is governed by the Riemann curvature, and that $\sigma^2=\id$ only when $M$ is Einstein.  We show that if $M$ has a conformal Killing vector field $\tau$ then  the cross product algebra $C(M)\rtimes_\tau\R$ viewed as a noncommutative analogue of $M\times\R$ has a natural $n+2$-dimensional calculus extending $\Omega^1$ and a natural spacetime Laplacian now directly defined by the extra dimension. The case $M=\R^3$ recovers the Majid-Ruegg bicrossproduct flat spacetime model and the wave-operator used in its variable speed of light preduction, but now as an example of a general  construction.  As an application we construct the wave operator on a noncommutative  Schwarzschild black hole and take a first look at its features. It appears that the infinite classical redshift/time dilation factor at the event horizon is made finite.  \end{abstract}
\maketitle 

\section{Introduction} 

This paper is divided into two parts.  Chapters 2,3 are mathematical and introduce a class of noncommutative or `quantum' spacetimes that are versions of $M\times \R$, where `space' is an undeformed classical Riemannian manifold $(M,\bar g)$ equipped with a conformal Killing vector field $\tau$. The `coordinate algebra' here is a cross product $C(M)\rtimes_\tau\R$ of the functions $C(M)$ on $M$ by the action of the vector field, i.e.  only noncommutative in respect of the time variable $t$ that governs the $\R$ direction. We also have as input data an arbitrary invertible function $\beta\in C(M)$ corresponding to a classical metric on $M\times\R$ of the static form 
\[ \bar g_{M\times \R}=\beta^{-1}\bd t\tens\bd t+\bar g.\]
Our main result, Corollary~3.4, is the construction of a natural noncommutative spacetime wave operator $\square$ on $C(M)\rtimes_\tau\R$ as a deformation of the classical spacetime Laplace-Beltrami operator for $\bar g_{M\times\R}$. This requires (Theorem~3.1) the construction of a locally inner quantum differential calculus on $C(M)\rtimes_\tau\R$ built in terms of the classical differential geometry on $M$ and finite difference geometry on $\R$. Our efforts go up to and including construction of the wave operator, without more non-essential further development of the noncommutative geometry. 

Our approach to the construction of the wave operator is a novel one which we believe to be deeper than mere analogy with the classical case. In fact, our point of view is that the conventional picture of physics emerges as a classical limit of a purely quantum/algebraic phenomena and not the other way around. Specifically, the origin of our approach is the discovery of a quantum anomaly for differential structures \cite{BegMa1,BegMa2,Ma:spo}  in the quantum groups approach to noncommutative geometry;  it was found that a range of sufficiently noncommutative geometries do {\em not} admit a covariant quantum differential calculus of classical dimensions but rather require at least one  `extra' cotangent direction $\theta$. The partial derivatives associated with such a $\theta$ turned out  to be noncommutative versions of the relevant Laplace-Beltrami operator. The phenomenon is not limited to quantum groups because most sufficiently noncommutative calculi can be expected to be inner in the sense of the existence of a 1-form $\theta$ such that $[\theta,f]=\lambda\extd f$ for any $f$ in the noncommutative `coordinate algebra' (or more generally in the exterior algebra and with a graded commutator). Here $\lambda$ is the noncommutativity parameter and both sides of this equation are zero as $\lambda\to 0$. Our approach is basically to turn this around and to define the Laplace operator on a noncommutative space as the partial derivative associated to such a $\theta$.  It is not known if such a  `spontaneous evolution'  is related to the modular automorphism group  which provides a canonical evolution of a von Neumann algebra more applicable in Connes' approach to noncommutative geometry \cite{Con}.

The present models are a little different due to being not very noncommutative. The main difference is a separation of the roles of the generator of the calculus, $\theta$, and the conjugate 1-form to the Laplacian, which we call $\theta'$. Our strategy is then as follows: starting with any classical Riemannian manifold $(M,\bar g)$ as `space' we use an operator $\bDelta$ (which can be taken to be the Laplace-Beltrami on $(M,\bar g)$) to define a natural extra-dimensional noncommutative differentuial calculus  such that $\extd f = \bd f + {\lambda\over 2}(\bDelta f)\theta' $. We then extend the calculus to one on $C(M)\rtimes_\tau \R$ (Theorem~3.1) by $\extd t$ and we then apply our principle backwards by defining the spacetime noncommutative Laplace-Beltrami $\square$ by 
\[ \extd \psi=\bd \psi+(\del_0\psi)\extd t+{\lambda\over 2}(\square\psi)\theta'\]
for all $\psi\in C(M)\rtimes_\tau\R$ (Proposition~3.3 and Corollary~3.4). The spacetime calculus is still inner locally (Proposition~3.6) and the difference between its generator $\theta$ and $\extd t$ is  proportional to $\theta'$.  We do not think of $\theta'$ as directly related any more to time but there are indications\cite{FreMa} that it may be related to the renormalisation group flow in quantum gravity.

The preliminary Section~2 undertakes the development of the spatial noncommutative calculus on $(M,\bar g)$ mentioned above, i.e. before the introduction of the time variable. It should be of interest in its own right as a natural way of encoding the structures of classical Riemannian geometry  into a natural `almost commutative' geometry. As the coordinate algebra remains $C(M)$ the  noncommutativity enters only at the level of noncommutation of 1-forms with functions. The differential calculus we use on $C(M)$ seems first to have been noted in \cite{MH}, as well as subsequently connected to Ito stochastic differentials. We provide (Proposition~2.4) a natural bimodule connection $\nabla$ on the noncommutative $\Omega^1(M)$.  The main idea of a bimodule connection is that $\nabla$ is a left-connection (in the usual sense of the derivation rule) on the 1-forms $\Omega^1$ but also $\sigma(\omega,\extd f):=\nabla(\omega f)-(\nabla\omega)f$ should be well-defined for any $f$ in the `coordinate algebra' and any 1-form $\omega$. The notion first appeared in works \cite{Mou,DV1,DV2}; see \cite{BegMa4} for some recent work and further references. Our main result for geometers is (Proposition~2.7)  that the `generalised braiding' $\sigma$ for our canonical bimodule connection  encodes the classical Riemann curvature of $M$ as failure of the `quantum Yang Baxter' or braid relations, and encodes the Ricci tensor as failure of the generalised braiding to be involutive. This provides a novel point of view on classical Riemannian geometry suggesting an associated braided 2-category. It also allows us to define a natural $\Omega^2$ in this context which we need in Section~3. Its extension to all degrees should provide a noncommutative geometry point of view on classical Hodge theory. 
 
The second part of the paper comprising Sections 4,5 is physical model building using our framework, up to and including the Schwarzschild black hole. The preliminary Section 4.1 reformulates the classical theory normally done with polar/angular coordinates but now using algebraic projective module methods that do not refer to trigonometry. We then are then able in Section~4.2 to apply the machinery of Section~3 to quantise any spherically symmetric static metric including, in particular, a `black hole differential algebra' associated to a natural conformal Killing vector on the classical black hole 3-geometry.  

Our full black hole differential algebra is, however, rather complicated and further computations are done in Section~5 in a simplified model. Here the initial geometry is $M=\R^3$ (or rather with the origin deleted) and its quantisation for the conformal killing vector $\tau=r{\del\over\del r}$ recovers the standard `bicrossproduct model' quantum spacetime\cite{MaRue} previously coming out of quantum group theory but now arising more geometrically. Our approach to its differential calculus also has new functional parameter $\beta$ which we view as encoding a gravitational potential. Section 5.2 then uses a `minimal coupling' process to amplify from the flat 3-space geometry to the one for the Schwarzschild or any other static metric of the form $\beta^{-1}\bd t\tens\bd t+\bar g$.  Looking in detail at the Schwarzschild case we study the resulting noncommutative wave operator and argue  that the infinite redshift or time dilation factor which classically applies to photons emitted just outside the black hole event horizon is rendered finite in our noncommutative model (see the discussion below Proposition~5.3), i.e the classical coordinate singularity at the even horizon is smoothed out in the noncommutative geometry. This is then confirmed by numerical solutions (Figure 1) which indeed indicate a finite redshift at the horizon. We also find an `interregnum' region between the  exterior of the horizon and the boundary of the black hole interior (the two sides of the classical horizon now become separated), and we take a numerical look at  `standing waves' in the black hole interior. Although the present paper is mainly mathematical, these initial results provide a first impression of some of the physical implications of the model. 

Both of our black hole versions have little in common with a previous attempt to define quantum black holes\cite{Schupp} by Drinfeld-type twist, in which each of the spheres of fixed radius is noncommutative. This is orthogonal to our constructions (in our case $M$ is undeformed) and also contains\cite{BegMa3}  a hidden nonassociativity in the differential geometry even if the coordinate algebra happens to be associative.   Moreover,  the bicrossproduct model in 2+1 dimensions arises naturally in a certain limit of $2+1$ quantum gravity and the latter may similarly yield noncommutative versions $C(M)\rtimes_\tau\R$  where $M$ is a Riemann surface, thus going beyond the local `model spacetime' picture overviewed in \cite{MaSch}. Our new framework  may also provide a setting to explore certain ideas for  Einstein's equations beyond the quantum group toy models proposed in \cite{Ma:pla}. 

In Sections~2,3 we work over $C(M)$ taken as, say, the smooth functions on a smooth manifold $M$. For the time variable $t$ we work mainly with polynomials and commutation relations (i.e. as operators on the classical function and 1-form spaces) but we suppose that our constructions extend to other classes of functions of $t$ including $e^{\imath\omega t}$ when we come to  physical applications in later sections. 

\section{Almost commutative cotangent bundles}

Let $M$ be a Riemannian manifold with coordinate algebra $C(M)$ and classical exterior algebra $\bar{\Omega}, \bar\extd$ (we use bar to denote the classical as we shall shortly introduce a new one). We let $\(\ ,\ \)$ be the inverse of the classical metric $\bar g$ and $\bar\Delta$ a classical operator such that the polarization formula 
\begin{equation}\label{bDeltapol}\bDelta(fg)=(\bDelta f)g+f(\bDelta g)+2\(\bd f,\bd g\)
\end{equation}
holds for all $f,g\in C(M)$ (and later a similar operator on 1-forms). One can take $\bDelta$ throughout to be the classical Laplace-Beltrami operator but we will need a little more generality for our application in Section~3. The following construction appears  first to have been noted in \cite{MH} (in the Laplace-Beltrami case) and is also related to stochastic calculus  on a Riemannian manifold. In order to be self-contained and to have it in the generality we need, we include a direct proof. We recall that in noncommutative geometry a differential structure can be defined algebraically as a bimodule $\Omega^1$ of `1-forms' over the coordinate algebra and a map $\extd$ from the latter to $\Omega^1$ obeying the Leibniz rule, see \cite{Con}. This is more general than the classical notion of differential structure even when the algebra is commutative. 

\begin{lemma} Let $M$ be a Riemannian manifold with notations as above. Then $\Omega^1=\bar\Omega^1\oplus C(M)\theta'$ with $\theta'$ central and 
\[ f\bullet\omega=f\omega,\quad \omega\bullet f=\omega f+{\lambda}\(\omega,\bar\extd f\)\theta',\quad  \extd f= \bar\extd f + {\lambda\over 2} (\bar\Delta f)\theta'\]
for all $\omega\in\bar\Omega^1, f\in C(M)$ makes $(\Omega^1,\bullet,\extd)$ a noncommutative first order differential calculus over $C(M)$. The bimodule structure enjoys commutation relations
\[ [\omega,f]=\lambda\(\omega,\bar\extd f\)\theta',\quad [\theta',f]=0\]
where the new product is understood.
\end{lemma}
\proof We check that the algebra $C(M)$ acts from each side. Thus
\begin{eqnarray*}( \omega\bullet f)\bullet g&=& (\omega f)\bullet g+{\lambda}\(\omega,\bar\extd f\)\theta' g=\omega f g+{\lambda}\(\omega f,\bar\extd g\)\theta'+{\lambda}\(\omega,\bar\extd f\)g\theta'\\
\omega\bullet(fg)&=&\omega fg+{\lambda}\(\omega,\bar\extd (fg)\)\theta'=\omega fg+{\lambda}\(\omega,(\bar\extd f)g+f\bar\extd g\)\theta'\end{eqnarray*}
using the Leibniz rule for $\bar\extd$. The two expressions are equal by tensoriality of $\(\ ,\ \)$ allowing us to move $f$ and $g$ out.  We have to verify that we have a bimodule
\begin{eqnarray*}(f\bullet \omega)\bullet g&=&(f\omega)\bullet g=f\omega g+{\lambda}\(f\omega,\bar\extd g\)\theta'\\
f\bullet(\omega\bullet g)&=&f\bullet(\omega g+{\lambda}\(\omega,\bar\extd g\)\theta')=f\omega g+{\lambda}f\(\omega,\bar\extd g\)\theta'
\end{eqnarray*}
which are again equal by tensoriality. Finally, we verify that $\extd$ is a derivation:
\begin{eqnarray*} \extd(fg)&=&\bar\extd(fg)+{\lambda\over 2}\bar\Delta(fg)\theta'=(\bar\extd f)g+f\bar\extd g+{\lambda\over 2}((\bar\Delta f)g+f\bar\Delta g)\theta'+\lambda\(\bar\extd f,\bar\extd g\)\theta'\\
&=&(\bar\extd f)g+f\bar\extd g+{\lambda\over 2}((\bar\Delta f)\theta' g+f\bar\Delta g\theta')+\lambda\(\bar\extd f,\bar\extd g\)\theta'=\extd f\bullet g+f\bullet\extd g\end{eqnarray*}
from the definitions. Note that the product on the free bimodule spanned by central element $\theta'$ is that of $C(M)$ and us not deformed in the construction. We used a polarisation property of $\bar\Delta$ (which can easily be proven in local coordinates in the case of the Laplace-Beltrami operator from symmetry of the metric tensor used in defining the 2nd order differential operator). Note that one normally also requires $f\tens g\to f\extd g$ to be surjective and this may require further conditions on the Riemannian manifold.  
 \endproof

 One can also set up the bimodule symmetrically with a $\lambda/2$ modification from either side. Next we recall that a linear connection can be defined in noncommutative geometry abstractly as a map $\nabla:\Omega^1\to \Omega^1\hat{\tens}\Omega^1$ such that $\nabla(f\omega)=\extd f\hat{\tens}\omega+f\nabla\omega$ for all $f$ in our coordinate algebra and all  1-forms $\omega$. Here we use hats to stress that the tensor product is with respect to the bimodule structure, but we will omit the hats when the context is clear. As $\Omega^1$ is a bimodule we may have the luxury of an additional derivation property from the other side,
 \[ \nabla(\omega\bullet f)=(\nabla\omega)\bullet f+\sigma(\omega\hat{\tens}\extd f),\quad \sigma:\Omega^1\hat{\tens}\Omega^1\to \Omega^1{\hat{\tens}} \Omega^1.\]
If $\sigma$ exists it will be defined by this equation and will be  a bimodule map, but in general it need not exist. The definition goes back to \cite{Mou,DV1,DV2} and several subsequent works and constitutes an alternative to the quantum group frame bundles approach to Riemanian geometry in \cite{Ma:rie,Ma:rieq}. We refer to \cite{BegMa4} for a full set of references to the literature. The following lemma is the key to all that follows. Note that we can use the inverse metric to convert a classical 1-form $\omega$ to a vector field $\omega^*=\(\omega,\ )$ and pull this back to extend the action of the classical Levi-Civita covariant derivative to 1-forms by  $\bnabla_\omega:=\bnabla_{\omega^*}$.

\begin{lemma} \label{keylem} There is a well-defined left module map
\[  \phi:\bar\Omega^1\bar\tens\bar\Omega^1\to \Omega^1\hat\tens\Omega^1,\quad  \phi(\omega\bar\tens\eta)=\omega\hat{\tens}\eta -\lambda \theta'\htens\bnabla_{\omega}\eta,\quad\forall\omega,\eta\in\bar\Omega^1.\]
 Here $\bar\nabla$ is the classical Levi-Civita connection. Moreover,
 \[ \phi(\omega\bar\tens\eta)\bullet f=f\phi(\omega\bar\tens\eta)+\lambda\(\eta,\bd f\)\omega\htens\theta'+\lambda\theta'\htens\(\omega,\bd f)\eta+\lambda^2\(\bar\nabla\bd f)(\omega,\eta)\theta'\htens\theta'\]
 where we evaluate against the two outputs of $\bar\nabla$ using the inverse metric. We use the same formula to define $\phi$ inductively as a map $\bar\Omega^1{}^{\bar\tens n}\to \Omega^1{}^{\hat\tens n}$ for $n\ge 1$ by
 \[ \phi(\omega\bar\tens\eta\bar\tens\zeta\bar\tens\cdots)=\omega\hat\tens\phi(\eta\bar\tens\zeta\bar\tens\cdots)-\lambda\theta'\hat\tens\phi(\bar\nabla_\omega(\eta\bar\tens\zeta\bar\tens\cdots))\]
 \end{lemma}
 \proof We prove the inductive version. The $n=1$ case is our existing identification $\bar\Omega^1\subset\Omega^1$ as a left module and a different right module structure. In general
 \begin{eqnarray*} \phi(\omega f\bar\tens\eta\bar\tens\zeta\cdots)&=& \omega f\hat\tens\phi(\eta\bar\tens\zeta\cdots)-\lambda\theta'\hat\tens\phi(\bar\nabla_{f\omega}(\eta\bar\tens\zeta\cdots))\\
 &&= \omega \hat\tens f\phi(\eta\bar\tens\zeta\cdots)-\lambda\(\omega,\bd f\)\theta' \hat\tens \phi(\eta\bar\tens\zeta\cdots)-\lambda\theta'\hat\tens\phi(f\bar\nabla_{\omega}(\eta\bar\tens\zeta\cdots))\\
 &&= \omega \hat\tens \phi(f \eta\bar\tens\zeta\cdots)-\lambda\theta'\hat\tens\phi(\(\omega,\bd f\)\eta\bar\tens\zeta\cdots)+f\bar\nabla_{\omega}(\eta\bar\tens\zeta\cdots))\\
 &&=\phi(\omega\bar\tens f(\eta\bar\tens\zeta\cdots))
 \end{eqnarray*}
 assuming that $\phi$ is a well-defined left module map on $\eta\bar\tens\zeta\bar\tens\cdots$. It is clearly then also a left module map as $f\phi(\omega\bar\tens\eta\bar\tens\zeta\cdots)=f\omega\hat\tens\phi(\eta\bar\tens\zeta\bar\tens\cdots)-\lambda\theta'\hat\tens f\phi(\bar\nabla_\omega(\eta\bar\tens\zeta\bar\tens\cdots))=\phi(f\omega\bar\tens\eta\bar\tens\zeta\cdots)$. Finally, we compute the right module structure for $n=2$ (the general case is similar),
 \begin{eqnarray*}\phi(\omega\bar\tens \eta)\bullet f&=&\omega\hat\tens\eta\bullet f-\lambda\theta'\hat\tens(\bnabla_\omega\eta)\bullet f\\
 &=&\omega\hat\tens\eta f+\lambda\omega\htens\(\eta,\bd f)\theta'-\lambda\theta'\hat\tens(\bnabla_\omega\eta) f-\lambda^2\(\bd f,\bnabla_\omega\eta\)\theta'\htens\theta'\\
 &=&\phi(\omega\hat\tens\eta f\)+\lambda\theta'\htens\(\omega,\bd f\)\eta+\lambda\omega\bullet \(\eta,\bd f\)\htens\theta'-\lambda^2\(\bd f,\bnabla_\omega\eta\)\theta'\htens\theta'\\
 &=&\phi(\omega\hat\tens\eta f\)+\lambda\theta'\htens\(\omega,\bd f\)\eta+\lambda\omega \(\eta,\bd f\)\htens\theta'+\lambda^2(\(\omega,\bd \(\eta,\bd f\)\)-\(\bd f,\bnabla_\omega\eta\))\theta'\htens\theta'\\
 &=&\phi(\omega\hat\tens\eta f\)+\lambda\theta'\htens\(\omega,\bd f\)\eta+\lambda\omega \(\eta,\bd f\)\htens\theta'+\lambda^2\(\bnabla_\omega\bd f,\eta\)\theta'\htens\theta'
 \end{eqnarray*}
 using in the last step that the metric is compatible with $\bnabla$. 
 \endproof
 
 In the sequel we will frequently view a classical tensor in $K\in \bar\Omega^1\bar\tens\bar\Omega^1$ as a tensorial (i.e., module) map on 1-forms by $K\omega=(\id\tens\(\omega,\ \))K$ and as a tensorial bilinear on 1-forms by $K(\omega,\eta)=\(\omega,K\eta\)$ (in other words by `raising indices'). We will denote the transpose by $K^T$.  We will also need an extension of $\bDelta$ to 1-forms such that
 \begin{equation}\label{bDeltapol1}\bDelta(f\omega)=(\bDelta f)\omega+f\bDelta\omega+2\bnabla_{\bd f}\omega\end{equation}
 \begin{equation}\label{bDeltapol2}\bDelta(\(\omega,\eta)\)=\(\bDelta\omega,\eta\)+\(\omega,\bDelta\eta\)+2\(\bnabla\omega,\bnabla\eta\)\end{equation}
 \begin{equation}\label{ricciD}[\bDelta,\bd]f={\rm Ricci_{\bDelta}}(\bd f)\end{equation}
 for all $f\in C(M)$ and $\omega,\eta\in \bar\Omega^1(M)$ and some tensorial operator which we have denoted ${\rm Ricci}_{\bDelta}$. Here the inverse metric is extended to tensor products in the obvious way.  One can take here $\bDelta$ the Laplace-Beltrami operarator for which the three identities are easily proven in local coordinates and ${\rm Ricci}_\bDelta$ is the usual Ricci tensor. In this case the third identity is also clear if one notes that $\bDelta-{\rm Ricci}$ then coincides with the Hodge Laplacian (as an example of a Weitzenbock identity), and this commutes with $\bd$. In keeping with our emphasis on wave operators in this paper, one could regard this apparently little-known identity (\ref{ricciD}) as a definition of ${\rm Ricci}$ in a manner that brings out its physical significance.  
  
 \begin{lemma}\label{nabla} For any classical tensor $K:\bar\Omega^1\to\bar\Omega^1$ the classical Levi-Civita connection induces a left connection on $\Omega^1$ with
\[ \nabla\omega=\phi(\bnabla\omega)+{\lambda\over 2}\theta'\htens(\bar\Delta-K)\omega,\quad\forall\omega\in \bar\Omega^1\subset\Omega^1.\] 
\end{lemma}
\proof  Using Lemma~\ref{keylem} we have
\begin{eqnarray*} \nabla(f\omega)&=&\phi(\bnabla(f\omega))+{\lambda\over 2}\theta'\htens(\bar\Delta-K)(f\omega)\\
&=&\phi(\bd f\btens\omega+f\bnabla\omega)+{\lambda\over 2}\theta'\htens((\bDelta f)\omega+f(\bDelta-K)\omega+2\bnabla_{\bd f}\omega)\\
&=&f\nabla\omega+\bd f\htens\omega+{\lambda\over 2}\theta'\htens(\bDelta f)\omega=f\nabla\omega+\extd f\htens\omega
\end{eqnarray*}
Note that explicitly,
\[ \nabla\omega=\bar\nabla_1\omega\hat\tens\bar\nabla_2\omega-\lambda\theta'\hat\tens\left(\bar\nabla_{\bar\nabla_1\omega}\bar\nabla_2\omega - {1\over 2}(\bar\Delta-K)\omega\right)\]
where $\bar\nabla_1\omega\tens\bar\nabla_2\omega$ denotes a lift of $\bar\nabla$ from $\bar\tens$ to the vector space tensor product $\tens$, and we project this down to $\htens$. However, we shall endeavour to avoid such expressions by working via the properties of $\phi$. Note also that the value of $\nabla\theta'$ is left unspecified but we will be led to some natural choices for it later on. Then we define $\nabla(f\theta')=\extd f\htens\theta'+f\nabla\theta'$. 

\begin{proposition}\label{sigma} Suppose that $[\nabla\theta',f]=0$ for all functions $f$. Then
\[ \sigma(\omega\hat\tens\eta)=\eta\hat\tens\omega+\lambda\bar\nabla_\omega\eta\hat\tens\theta'-\lambda\theta'\hat\tens\bar\nabla_\eta\omega+\lambda(\omega,\eta)\nabla\theta'+{\lambda^2\over 2}({\rm Ricci}_\bDelta+K^T)(\omega,\eta)\theta'\hat\tens\theta'\]
\[ \sigma(\theta'\hat\tens\omega)=\omega\hat\tens\theta',\quad  \sigma(\omega\hat\tens\theta')=\theta'\hat\tens\omega,\quad  \sigma(\theta'\hat\tens\theta')=\theta'\hat\tens\theta'\]
for all $\omega,\eta\in\bar\Omega^1$ makes $\nabla$ in Lemma~\ref{nabla} into a bimodule connection.
\end{proposition}
\proof Note that we can write
\begin{equation}\label{sigmaphi}  \sigma(\omega\hat\tens\eta)=\phi(\eta\btens\omega)+\lambda\bar\nabla_\omega\eta\hat\tens\theta'+\lambda(\omega,\eta)\nabla\theta'+{\lambda^2\over 2}({\rm Ricci}_\bDelta+K^T)(\omega,\eta)\theta'\hat\tens\theta'.\end{equation}
We first check that this is well-defined as a left module map. In fact this works for any tensor $R$ in the role of ${\rm Ricci}_\bDelta+K^T$,
\begin{eqnarray*}
\sigma(\omega \htens f\eta)&=&\phi(f\eta\btens\omega)+\lambda f\bnabla_\omega\eta+\lambda\(\omega,\bd f\)\eta\htens\theta'+
f\lambda\(\omega,\eta\)\nabla\theta'+f{\lambda^2\over 2}R(\omega,\eta)\\
&=&f\sigma(\omega\htens \eta)+\lambda(\omega,\bd f)\eta\htens\theta'=\sigma(\omega\bullet f\hat\tens\eta)
\end{eqnarray*}
provided we define $\sigma(\theta'\htens\eta)=\eta\htens\theta'$. This indeed provides the bimodule connection property on $\nabla\theta'$ provided $[\nabla\theta',f]=0$ and provided we define $\sigma(\theta'\hat\tens\theta')=\theta'\htens\theta'$. As $\theta'$ commutes with functions, our assumption on $\nabla\theta'$ ensures that everything in this sector behaves as classically. Finally, we check the bimodule connection property. Using the definitions and the commutation relations as well as Lemma~\ref{keylem},
\begin{eqnarray*}
\nabla(\omega\bullet f)-(\nabla\omega)\bullet f&=&\nabla(\omega f)+\lambda\extd\(\bd f,\omega\)\htens\theta'+\lambda\(\bd f,\omega\)\nabla\theta'-\phi(\bnabla\omega)\bullet f-{\lambda\over 2}\theta'\htens((\bDelta-K)\omega)\bullet f\\
&=&\nabla(\omega f)+\lambda\bd\(\bd f,\omega\)\htens\theta'+{\lambda^2\over 2}\bDelta\(\bd f,\omega\)\htens\theta'
+\lambda\(\bd f,\omega\)\nabla\theta'\\
&&-\phi((\bnabla\omega) f)-\lambda\(\bd f,\bnabla\omega)\htens\theta'-\lambda\theta'\htens\bnabla_{\bd f}\omega-{\lambda^2}(\bnabla\bd f)(\bnabla\omega)\\
&&-{\lambda\over 2}\theta'\htens((\bDelta-K)\omega)f-{\lambda^2\over 2}\(\bd f,(\bDelta-K)\omega\)\\
&=& \phi(\bnabla(\omega f)-(\bnabla\omega)f)+{\lambda\over 2}\theta'\htens\left((\bDelta-K)(\omega f)-((\bDelta-K)\omega)f-2\bnabla_{\bd f}\omega\right)\\
&&+\lambda\(\bnabla\bd f,\omega\)\htens\theta'+\lambda\(\bd f,\omega\)\nabla\theta'+{\lambda^2\over 2}\(\bd f, K\omega\)\\
&&+{\lambda^2\over 2}\left(\bDelta\(\bd f,\omega\)-2(\bnabla\bd f)(\bnabla \omega)-\(\bd f,\bDelta\omega\)   \right)\theta'\htens\theta'\\
&=&\phi(\bd f\btens\omega)+{\lambda\over 2}\theta'\htens(\bDelta f)\omega+\lambda\bnabla_\omega\bd f\htens\theta' +\lambda\(\bd f,\omega\)\nabla\theta'+{\lambda^2\over 2}\(\bd f, K\omega\)\\
&&+{\lambda^2\over 2}\(\bDelta\bd f,\omega)\theta'\htens\theta'\\
&=&\sigma(\omega\htens\bd f)+{\lambda\over 2}\theta'\htens(\bDelta f)\omega+{\lambda^2\over 2}\(\bd\bDelta f,\omega)\theta'\htens\theta'\\
&=&\sigma(\omega\htens\bd f)+{\lambda\over 2}\theta'\htens\omega\bullet\bDelta f=\sigma(\omega\htens\extd f)
\end{eqnarray*}
provided $\sigma(\omega\htens\theta')=\theta'\htens\omega$, because in this case have a right module property
$ \sigma(\omega\htens\theta' f)=\sigma(\omega\bullet f\htens\theta')=\sigma(f\omega\htens\theta')+\lambda\(\omega,\bd f\)\theta'\htens\theta'=f\theta'\htens\omega+\lambda\(\omega,\bd f\)\theta'\htens\theta'=\theta'\htens\omega\bullet f=\sigma(\omega\htens\theta')\bullet f$ for any function $f$ (here applied to $\bDelta f$). 
In the main computation we used that $\bnabla$ is metric compatible to compute $\bd\(\bd f,\omega)$ and torsion free in the form $\(\bnabla \bd f,\omega\)= \bnabla_\omega\bd f$ (this is because antisymetrization of $\bnabla\omega$ is provided by $\wedge$ and $\bnabla\wedge\bd f=\bd \bd f=0$; we will explain this point of view on torsion further below) and we used   (\ref{ricciD}).  Although this only verifies  $\sigma(\omega\htens\extd f)$, the properties of $\nabla$ a connection imply that it is fully a right module map and that it is fully determined. Thus let  $\eta=a\bd b$ (or a sum of such terms) then $\sigma(\omega\htens a\bd b)=\sigma(\omega\bullet a\htens\bd b)=a\sigma(\omega\htens\bd b)+\lambda(\omega,\bd a)\bd b\htens\lambda\theta'=a\phi(\bd b\htens\omega)+\lambda\bnabla_\omega(a\bd b)+{\rm tensorial\ terms}=\sigma(\omega\htens\eta)$ in the form (\ref{sigmaphi}). \endproof

Let us note that $\sigma$ does not map over from the flip map under $\phi$, rather we can write the above result as
\begin{equation}\label{sigmaflip}\sigma(\phi(\omega\btens\eta))=\phi(\eta\btens\omega)+\lambda\(\omega,\eta)\nabla\theta'+{\lambda^2\over 2}({\rm Ricci}_\bDelta+K^T)(\omega,\eta)\theta'\tens\theta',\quad\forall\omega,\eta\in\bar\Omega^1.\end{equation}

We now consider the space of 2-forms. Any $(\Omega^1,\extd)$ on an algebra has a `maximal prolongation' obtained by the minimal requirements that $\extd$ extends as a graded derivation with $\extd^2=0$. This is typically too large except in almost-commutative cases. In our case we take the maximal prolongation modulo the relations
\begin{equation}\label{thetawedge} \{\omega,\theta'\}=\theta'^2=0,\quad\forall\omega\in \bar\Omega^1\end{equation}
to be consistent with our assumptions leading to the corresponding classical behaviour of $\sigma$. It remains to find the relations in $\Omega^2$ explicitly and $\extd$ on $\Omega^1$. This will be tied up with torsion and we recall that in terms of forms this can be written as\cite{Mou,Ma:rie},
\begin{equation}\label{torsion} T_\nabla(\omega):=\nabla\wedge\omega-\extd\omega,\quad T_\nabla:\Omega^1\to \Omega^2.\end{equation}
This is usually discussed in the context of a metric compatible connection but one can take it as a definition, both in the classical case where it applies to $\bnabla$, and in the `quantum case'.

\begin{proposition}\label{omega2} The relations $\wedge{\rm image}(\id+\sigma)=0$ for $\sigma$ in Proposition~\ref{sigma}, i.e. 
\[  \{\omega,\eta\}=\lambda\theta'(\bnabla_\omega\eta+\bnabla_\eta\omega)-\lambda\(\omega,\eta)\nabla\wedge\theta'\]
hold  for all $\omega,\eta$ in $\Omega^1$, provided $\nabla\wedge\theta'=\extd\theta'$. Moreover,
\[ \extd \omega=\nabla\wedge \omega\]
so that $\nabla$ in Proposition~\ref{nabla} has zero torsion, provided  $\bDelta$ is the Laplace-Beltrami operator and $K={\rm Ricci}$.
\end{proposition}
\proof We apply $\extd$ to the relations in degree 1 under the assumption that $\extd^2=0$ to obtain
\[ \{\theta',\extd f\}=[\extd\theta',f],\quad \{\extd b,\extd f\}+\lambda(\extd\(\bd f,\bd b\))\theta'+\lambda\(\bd b,\bd f\)\extd\theta'=0.\] 
We note in passing that the first equation means $\{\theta',a_i\bd b_i\}=a_i[\extd\theta',b_i]-\lambda a_i(\bDelta b_i)\theta'^2$ which will not depend only on $\omega=a_i\extd b_i$ unless the right hand side is zero, which in turn implies that  
\begin{equation}\label{thetaf} [\extd\theta',f]=0,\quad \forall f\in C(M)\end{equation}
and hence that $\theta'^2=0$. Hence (\ref{thetawedge}) are the only reasonable assumptions for the calculus to be `built on' the classical one. In the second equation, assuming (\ref{thetawedge}), we can replace $\extd$ by $\bd$. Then
\begin{eqnarray*}\{a_i\bd b_i,\bd f\}&=&a_i\{\bd b_i,\bd f\}+[\bd f,a_i]\bd b_i\\
&=&-\lambda a_i(\bd\(\bd b_i,\bd f\))\theta'-\lambda\(\omega,\bd f\)\extd\theta'+\lambda\(\bd a_i,\bd f\)\bd b_i\theta'\\
&=&-\lambda(\bd\(\omega,\bd f\))\theta'+\lambda\bd (a_i \(\bd b_i,\bd f\)-\(\bd a_i,\bd f\)\bd b_i)\theta'-\lambda\(\omega,\bd f\)\extd\theta'\\
&=&-\lambda((\id\bar\tens i_{\bd f})\bnabla\omega+(\id\tens i_{\omega})\bnabla\bd f)\theta'-\lambda(\bnabla_{\bd f}\omega-(\id\bar\tens i_{\bd f})\bnabla\omega)\theta'-\lambda\(\omega,\bd f\)\extd\theta'\\
&=&-\lambda\bnabla_\omega\bd f\theta'-\lambda\bnabla_{\bd f}\omega\theta'-\lambda\(\omega,\bd f\)\extd\theta'\\
\end{eqnarray*}
where we used metric compatibility to expand $\extd\(\omega,\bd f\)$ and recognised the interior product $i_{\bd f}\bd \omega=i_{\bd f}\bnabla\wedge\omega$ using zero torsion (in the form (\ref{torsion})). Here $i$ denotes interior product. Finally, we used that $\bnabla\bd f$ is symmetric in its two outputs, again due to zero torsion and $\bd^2 f=0$. We then  let $\eta=a_i\bd b_i$ be some other 1-form and find
\begin{eqnarray*} \{\omega,a_i\bd b_i\}&=&[\omega,a_i]\bd b_i+a_[\{\omega,\bd b_i\}\\
&=&\lambda(\omega,\extd a_i)\theta'\bd b_i -\lambda a_i\bnabla_\omega\bd b_i\theta'-\lambda\bnabla_{\eta}\omega\theta'-\lambda\(\omega,\eta\)\extd\theta'\\
\end{eqnarray*}
which gives the expression stated provided $\nabla\wedge\theta'=\extd\theta'$, i.e. that $T_\nabla(\theta')=0$. This completes the first part of the proposition. 

For the second part, it is sufficient to prove that $\nabla\wedge\extd f=0$ for all $f$ as then $\nabla(a_i\extd b_i)=\extd a_i\wedge\extd b_i=\extd(a_i\extd b_i)$. Thus,
\begin{eqnarray*} \nabla\wedge \extd f&=&\nabla\wedge\bd f+{\lambda\over 2}(\bd\bDelta f\, \theta'+\bDelta f\nabla\wedge\theta')\\
&=&\nabla\wedge \bd f+{\lambda\over 2}((\bDelta-{\rm Ricci}_\bDelta)\bd f\, \theta'+\bDelta f\nabla\wedge\theta')\\
&=&\wedge\phi(\bnabla\bd f)+{\lambda\over 2}(\theta' (\bDelta-K)\bd f+(\bDelta-{\rm Ricci}_\bDelta)\bd f\, \theta'+\bDelta f\nabla\wedge\theta')\\
&=&\wedge\phi(\bnabla\bd f)+{\lambda\over 2}\bDelta f\nabla\wedge\theta' \end{eqnarray*}
provided $K={\rm Ricci}_\bDelta$. We used (\ref{ricciD}) and the definitions. Assuming this, and proceeding further we compute from (\ref{sigmaflip}),
\begin{eqnarray*} \sigma(\phi(\bnabla\bd f))&=&\phi(\bnabla\bd f)+\lambda(\(\ ,\ \)\bnabla\bd f) \nabla\theta'+{\lambda^2}{\rm Ricci}_\bDelta(\bnabla\bd f)\theta'\htens\theta'\\
(\id+\sigma)(\phi(\bnabla\bd f))&=&2\phi(\bnabla\bd f)+\lambda\bDelta_{LB} f \nabla\theta'+{\lambda^2}{\rm Ricci}_\bDelta(\bnabla\bd f)\theta'\htens\theta'
\end{eqnarray*}
since $\bnabla\bd f$ is symmetric in its two outputs (as the classical torsion vanishes) and  $\bDelta_{LB}=\(\ ,\ \)\bnabla\bd$ the Laplace-Beltrami. Applying $\wedge$ we deduce that $\wedge\phi(\bnabla\bd f)=-{\lambda\over 2}\bDelta_{LB} f\nabla\wedge\theta'$ and hence that $\nabla\wedge\extd f=0$ provided we take $\bDelta=\bDelta_{LB}$. 
\endproof

It is possible to take this as a definition, i.e using the Laplace-Beltrami operator and defining
\[ \Omega^2:=\Omega^1\htens\Omega^1/{\rm image}(\id+\sigma),\quad \extd=\nabla\wedge.\]
 It should also be possible to  proceed to construct an entire exterior algebra $(\Omega,\extd)$ using for $\extd$ the same format and the Hodge laplacian on forms. This will be considered elsewhere  as it is tangential to our main purpose (for Riemannian geometry we need mainly up to degree 2). 

We now turn to first consideration of the `quantum metric'. There is a natural `initial metric' $g$ which we consider but we note that $\nabla$ in our construction is only guaranteed to be metric compatible for it to $O(\lambda)$. Hence typically either $g$ or $\nabla$ will need to be modified to fully extend the Riemainnian geometry to the quantum case. Again, this is not essential to the noncommutative wave operator and we defer the topic to a sequel. 

\noindent \begin{corollary} Let
\[ g=\phi(\bar g)=\bar g_1\hat\tens\bar g_2-\lambda\theta'\hat\tens\bar\nabla_{\bar g_1}\bar g_2\]
where we take any lift $\bar g_1\tens\bar g_2$ of the classical metric $\bar g\in \Omega^1\bar\tens\Omega^1$.
\begin{enumerate} 
\item 
 $\sigma(g)=g$ iff $\nabla\theta'=-{\lambda\over 2}({\rm Ricci}_\bDelta+K^T)(\bar g)\theta'\htens\theta'$.
\item In the context of Proposition~\ref{omega2}, $\wedge(g)=0$ iff $\nabla\wedge\theta'=0$
\item $g$ generates the calculus in the sense 
\[ [g,f]=\lambda\extd f\htens\theta'+\lambda\theta'\htens\extd f,\quad\forall f\in C(M)\]
 provided that $\bDelta$ is the Laplace-Beltrami operator.
\end{enumerate}
\end{corollary}
\proof We have immediately from (\ref{sigmaflip}) and symmetry of the classical metric that 
\begin{eqnarray*}\sigma(g)
&=&g+\lambda\dim(M)\nabla\theta'+{\lambda^2\over 2}({\rm Ricci}_\bDelta+K^T)(\bar g)\theta'\htens\theta'
\end{eqnarray*}
 This gives the condition for $\sigma(g)=g$ and also implies from looking at $0=\wedge(\id+\sigma)(g)$ that 
 \[ \wedge(g)=-{\lambda\over 2}\dim(M)\nabla\wedge\theta'.\]
 Finally, the commutation relation follows immediately from Lemma~\ref{keylem} as this gives
\[ [g,f]=\lambda\(g_1,\bd f\)g_2\htens\theta+\lambda\theta'\htens\(g_1,\bd f\)g_2+\lambda^2(\bnabla\bd f)(\bar g)\theta'\htens\theta'\]
and in the last term we use $(\bnabla\bd f)(\bar g)=\bDelta f$, i.e. the  Laplace-Beltrami operator. \endproof

This suggests that $g$ is not yet the definitive noncommutative metric.  We conclude with some further properties of $\sigma$ not directly relevant to the our application but possibly of interest to classical geometry. We let  $\bDelta$ be the classical Laplace-Beltrami operator and $K={\rm Ricci}$.

\begin{proposition} Let $(M,\bar g)$ be a Riemannian manifold and $\sigma$ be the cannonical bimodule map $\sigma$ in Proposition~\ref{sigma} associated to the Laplace-Beltrami operator and $K={\rm Ricci}$. 
\begin{enumerate}
\item $\sigma^2=\id$ iff ${\rm Ricci}=\mu (\ ,\ \)$ (so $(M,\bar g)$ is Einstein) and ${1\over 2}(\id+\sigma)\nabla\theta'=-\lambda\mu\theta'\htens\theta'$ for some $\mu\in C(M)$.
\item Suppose that $\nabla\theta'=-\lambda\mu\theta'\htens\theta'$ for some $\mu\in C(M)$. Then $\sigma$ obeys the braid or `quantum Yang-Baxter' equations iff the $(M,\bar g)$ is flat and $\mu=0$.
\end{enumerate}
\end{proposition}
\proof Either from (\ref{sigmaflip}) or by direct computation we find that
\[\sigma^2(\omega\htens\eta)=\omega\htens\eta+\lambda\(\omega,\eta\)(\id+\sigma)\nabla\theta'+2\lambda^2{\rm Ricci}(\omega,\eta)\theta'\htens\theta',\quad\forall\omega,\eta\in\bar\Omega^1.\]
Hence $\sigma=\id$ obtains iff $ -\lambda{{\rm Ricci}(\omega,\eta)\over\(\omega,\eta\)}\theta'\htens\theta'={\id+\sigma\over 2}\nabla\theta' $
for all $\omega,\eta$ with  $\(\omega,\eta\)\ne 0$ (and ${\rm Ricci}(\omega,\eta)=0$ if $(\omega,\eta)=0$.) The right hand side does not depend on $\omega,\eta$ hence ${\rm Ricci}(\omega,\eta)=\mu\(\omega,\eta\)$ for some  $\mu\in C(M)$. For the second part, let $\sigma_{12}$ denote $\sigma$ acting on the first two factors of $\Omega^{\htens 3}$. We find after a lengthy but elementary computation, 
\begin{eqnarray*}(\sigma_{12}\sigma_{23}\sigma_{12}-\sigma_{23}\sigma_{12}\sigma_{23})(\omega\htens\eta\htens\zeta)&=& \lambda^2\theta'\htens\theta'\htens([\bnabla_\zeta,\bnabla_\eta]-\bnabla_{\bnabla_\eta\zeta-\bnabla_\zeta\eta})\omega\\
&&+ \lambda^2\theta'\htens([\bnabla_\omega,\bnabla_\zeta]-\bnabla_{\bnabla_\omega\zeta-\bnabla_\zeta\omega})\eta\htens\theta'\\
&&+ \lambda^2([\bnabla_\eta,\bnabla_\omega]-\bnabla_{\bnabla_\eta\omega-\bnabla_\omega\eta})\zeta\htens \theta'\htens\theta'+O(\lambda^3\theta'^{\htens 3})\end{eqnarray*}
so the QYBE hold to this order iff the full Riemann curvature vanishes. The omitted $\theta'^{\htens 3}$ terms involve the Ricci curvature on the one hand and terms involving $\mu$ on the other. The latter are, using metric compatibility, 
\begin{eqnarray*}&&\lambda^3\left(\(\bd\mu,(\eta,\zeta)\omega+(\omega,\eta)\zeta-(\omega,\zeta)\eta\)+2\mu(\(\omega,\bnabla_\zeta\eta\)+\(\zeta,\bnabla_\omega\eta\))\right)\theta'^{\htens 3}
\end{eqnarray*}
This cannot vanish for all $\omega,\eta,\zeta$ unless $\mu=0$. For example, set $\omega=\eta=\zeta$ and $\eta$ such that $\bnabla_\eta\eta=0$ at a point $x\in M$ and with any chosen direction $\eta(x)$. Then the second term vanishes and we conclude that $\extd \mu=0$ at any point. We can then take $\eta$ such that $2(\eta,\bnabla_\eta\eta)=\eta(\extd(\eta,\eta))\ne0$ at any point to conclude that $\mu=0$ there. \endproof

The simplest case is to take $\nabla\theta'=0$. In this case $\sigma$ obeys the braid relations {\em iff} $(M,\bar g)$ is flat and is involutive {\em iff} it is Ricci flat. 

\section{Wave operator on $C(M)\rtimes\R$ as quantisation of $M\times \R$}

We are now going to use the machinery of the previous section to construct a noncommutative spacetime deforming $M\times \R$, a differential calculus and a wave operator $\square$ on it. As `coordinate algebra' we let $A=C(M)\rtimes\R$ where we adjoin a variable $t$ for `time', with relations 
\[ [f,t]=\lambda\tau(f)\]
where $\tau$ is a vector field on $M$. We have used the same deformation parameter as before but without loss of generality as we could change the normalisation of $\tau$. This algebra has a noncommutative time variable as with the bicrossproduct model spacetime and is manifestly associative because any vector field $\tau$ generates an infinitesimal action of $\R$ on the algebra $C(M)$ and our algebra is the semidirect product by this. At least when $M$ is compact one can exponentiate the action as well as complete to a $C^*$ algebra if one wishes, although we shall not do either of these steps here.  

In order to apply the theory of Section~2 we let $\bDelta_{LB}$ be the Laplace-Beltrami operator on $(M,\bar g)$ and $\zeta$ a classical vector field on $M$,  and define
\begin{equation}\label{bDeltazeta}\bDelta f= \bDelta_{LB} f+ \zeta(f),\quad \bDelta\omega=\bDelta_{LB} \omega+\bnabla_\zeta\omega\end{equation}
for all $f\in C(M)$ and $\omega\in \bar\Omega^1(M)$. One may check that the properties (\ref{bDeltapol}), (\ref{bDeltapol1}),  (\ref{bDeltapol2}), (\ref{ricciD}) continue to hold with
\[ {\rm Ricci}_\bDelta={\rm Ricci}+\bnabla_\zeta-\bar\CL_\zeta\]
as an operator on $\Omega^1$, where $\bar\CL_\zeta$ is the Lie derivative along $\zeta$. We will later fix $\zeta$ in terms of a functional parameter below, but for the moment it is unspecified. From Section~2 we have an extended differential calculus $(\Omega^1,\extd)$ and other structures constructed from $(M,\bar g,\zeta)$. We let $\zeta^*$ be the 1-form corresponding to $\zeta$ under the metric.

\begin{theorem} Let $M$ be a Riemannian manifold equipped with a vector field $\zeta$,  $\beta\in C(M)$ and  $\tau$ a conformal Killing vector field. Then the  calculus $(\Omega^1,\extd)$ on $M$ defined by $\zeta$ extends to a  first order differential calculus $(\Omega^1(C(M)\rtimes\R),\extd)$ with further relations
\[  [\omega,t]=\lambda(\bar\CL_\tau-\id) \omega-\lambda^2({\scriptstyle {n-2\over 4}})\(\bar\extd\alpha,\omega\)\theta'-{\lambda^2\over 2}\(\bar\CL_\tau\zeta^*,\omega\)\theta'\]
\[  [\theta',t]=\alpha\lambda \theta',\quad [f, \extd t]=\lambda \extd f,\quad [\extd t,t]=\beta\lambda\theta'-\lambda\extd t\]
for all $\omega\in\bar\Omega^1(M), f\in C(M)$. Here $n=\dim(M)$ and  $\alpha={2\over n}{\rm div}(\tau)-1$.
\end{theorem}
\proof That $\tau$ is a conformal Killing vector field can  be written in terms of the inverse metric as as
\begin{equation}\label{taueqn}  \tau(\(\omega,\eta\))=\(\bar \CL_\tau\omega,\eta\)+\(\omega,\bar \CL_\tau\eta\)-(1+\alpha)\(\omega,\eta\)\end{equation}
which is the form we shall use. We have to verify the various Jacobi identities concerning the extension by $t,\extd t$. Thus
\begin{eqnarray*}&&\kern-15pt \lambda^{-2}\left([[\omega,t],f]+[[t,f],\omega]+[[f,\omega],t]\right)\\
&&=\lambda^{-1}[\bar\CL_\tau(\omega)-\omega-\lambda({\scriptstyle {n-2\over 4}})\(\bar\extd\alpha,\omega\)\theta'-{\lambda\over 2}\(\bar\CL_\tau\zeta^*,\omega\)\theta',f]+\lambda^{-1} [\omega,\tau(f)]-\lambda^{-1}[\(\omega,\bar\extd f\)\theta',t]\\
&&=\(\bar\CL_\tau(\omega)-\omega,\bar\extd f\)\theta'+\(\omega,\bar\extd\tau(f)\)\theta'-\tau(\(\omega,\bar\extd f\))\theta'-\alpha\(\omega,\bar\extd f\)\theta'=0\end{eqnarray*}
as (\ref{taueqn}) with $\eta=\bar\extd f$. Also $[[\theta',t],f]+[[t,f],\theta']+[[f,\theta'],t]=0$ as each term is zero, while
\begin{eqnarray*}&&\kern-15pt [[\extd t,f],g]+[[f,g],\extd t]+[[g,\extd t],f]=-\lambda[\extd f,g]+\lambda[\extd g,f]=-\lambda[\bar\extd f,g]+\lambda[\bar\extd g,f]=0\end{eqnarray*}
by symmetry of $\(\ ,\ \)$. Finally,
\begin{eqnarray*}&&\kern-15pt \lambda^{-1}\left([[\extd t,t],f]+[[t,f],\extd t]+[[f,\extd t],t]\right)\\
&=&[\beta\theta'-\extd t,f]-[\tau(f),\extd t]+[\extd f,t]=\lambda\extd f-\lambda\extd\tau(f)+[\extd f,t]\\
&=&\lambda\extd f-\lambda\extd \tau(f)+[\bar\extd f+{\lambda\over 2}\bar\Delta f\theta',t]\\
&=&-\lambda(\extd \tau(f)-\extd f)+\lambda(\bar\extd\tau(f)-\bar\extd f) -\lambda^2({\scriptstyle{n-2\over 4}})\(\bar\extd \alpha,\bar\extd f\)\theta'-{\lambda^2\over 2}\(\bar\CL_\tau\zeta^*,\bd f\)\theta'\\
&&\quad+\alpha{\lambda^2\over 2}\bar\Delta f\theta'+{\lambda^2\over 2}\tau(\bar\Delta f)\theta'\\
&=&{\lambda^2\over 2}\left(-\bar\Delta\tau(f)+(1+\alpha)\bar\Delta f+\tau(\bar\Delta f)-({\scriptstyle{n-2\over 2}})\(\bar\extd \alpha,\bar\extd f\)-\(\bar\CL_\tau\zeta^*,\omega\)\right)\theta'=0
\end{eqnarray*}
by a property of conformal Killing vectors in the lemma below. Note that in this case we obtain a formula for the commutator with quantum differentials
\begin{equation}\label{dft} [\extd f,t]=\lambda(\extd \tau(f) -\extd f).\end{equation}
This is equivalent to checking that $\extd$ is a derivation for products with $t$ in the sense $\extd[f,t]=\lambda\extd \tau(f)=[\extd f,t]+[f,\extd t]$. Once again we do not verify surjectivity, however this appears to be true in practice. 

To complete the proof we need the following elementary lemma which extends a well-known property of Killing forms to the conformal case as well as to our slightly more general $\bDelta$. As we have not found it in the literature, we include a short proof for completeness.

\begin{lemma} If $M$ is a Riemannian manifold with Laplace-Beltrami operator $\bar\Delta_{LB}$ and $\tau$ a conformal Killing vector field, then
\[ [ \bar\Delta_{LB},\tau]f=(1+\alpha)\bar\Delta_{LB} f - {n-2\over 2}\(\bar\extd\alpha,\bar\extd f\)\]
for all $f\in C(M)$, where $1+\alpha= {2\over n}{\rm div}(\tau)$ in our previous conventions.
\end{lemma}
\proof We can do this using the abstract notation of Section~2 but since this is a classical result we will use standard local coordinate methods. In local form one can write the equations for a conformal Killing vector as $\bar\nabla_a\tau_b+\bar\nabla_b\tau_a=(1+\alpha)g_{ab}$ where $\tau_a=g_{ab}\tau^b$ and $\bar\nabla$ is the classical Levi-Civita connection. We assume that the vector fields of our local coordinate system commute. Then applying $\bar\nabla^a$ to both sides (summation understood) we have $\bar\Delta_{LB} \tau_b=-\bar\nabla^a\bar\nabla_b\tau_a+\bar\nabla_b\alpha=-R^a{}_b{}_a{}^c\tau_c-\bar\nabla_b{\rm div}(\tau)-\bar\nabla_b\alpha=-{\rm Ricci}_b{}^c\tau_c-({n-2\over 2})\bar\nabla_b\alpha$. Next we compute
\begin{eqnarray*} \bar\Delta_{LB}\tau f&=&g^{ab}\bar\nabla_a\bar\nabla_b(\tau_c\bar\nabla^c f)=g^{ab}\bar\nabla_a((\bar\nabla_b\tau_c)(\bar\nabla^c f))+g^{ab}\bar\nabla_a(\tau_c\bar\nabla_b\bar\nabla^c f)\\
&=&g^{ab}(\bar\nabla_a\bar\nabla_b\tau_c)\bar\nabla^cf+g^{ab}(\bar\nabla_b\tau_c)(\bar\nabla_a\bar\nabla^c f)+2g^{ab}(\bar\nabla_a\tau_c)(\bar\nabla_b\bar\nabla^c f)+g^{ab}\tau_c\bar\nabla_a\bar\nabla_b\bar\nabla^c f\\
&=&(\bar\Delta_{LB} \tau_c)\bar\nabla^cf+(1+\alpha)\bar\Delta_{LB} f+g^{ab}\tau_c\bar\nabla_a\bar\nabla^c \bar\nabla_b f\\
&=&(\bar\Delta_{LB} \tau_c)\bar\nabla^cf+(1+\alpha)\bar\Delta_{LB} f+g^{ab}\tau_c\bar\nabla^c\bar\nabla_a \bar\nabla_b f+g^{ab}\tau_c R_a{}^c{}_b{}^d\bar\nabla_d f\\
&=&(\bar\Delta_{LB} \tau_c)\bar\nabla^cf+(1+\alpha)\bar\Delta_{LB} f+\tau\bar\Delta_{LB} f+\tau^c{\rm Ricci}^{cd}\bar\nabla_d f
\end{eqnarray*}
for all $f$. We use the Leibniz rule and that our local basis covariant derivatives commute when acting on functions. We then combine these two observations. 
\endproof

From the result for $\bar\Delta_{LB}$, the general case
\[ [ \bar\Delta,\tau]f=(1+\alpha)\bar\Delta f - {n-2\over 2}\(\bar\extd\alpha,\bar\extd f\)-\(\bar\CL_\tau\zeta^*,\bd f\)\]
then follows by an elementary computation and completes the proof of the theorem. 

The theorem provides a noncommutative geometry $\Omega^1(C(M)\rtimes\R)$ built on classical 1-forms with an extra cotangent direction $\theta'$ in addition to $\extd t$ and is a little more than an abstract calculus defined commutation relations among functions and {\em quantum} differentials. At the latter level the theorem looks simpler and we collect all the relations together for reference as
\begin{eqnarray}\label{cr} [f,g]=0,\quad&& [f,t]=\lambda\tau(f),\quad [\extd f,g]=\lambda\(\bar\extd f,\bar\extd g\)\theta',\quad [\theta',f]=0,\quad [\theta',t]=\alpha\lambda \theta'\nonumber \\ 
&& [\extd f,t]=\lambda(\extd\tau(f)-\extd f),\quad [f,\extd t]=\lambda\extd f,\quad [\extd t,t]=\beta\lambda\theta'-\lambda\extd t.\end{eqnarray}
for all $g,f\in C(M)$.    Having obtained this structure, one could take these relations as a definition of the calculus and verify the Jacobi identities, one of which would rapidly  lead back to the conformal Killing equation (\ref{taueqn}). Our more involved proof of Theorem~3.1 shows that the construction is properly defined with respect to the structure of the manifold $M$ by virtue of  being built on the classical objects and it is only there that the choice of $\zeta$ is visible.

The case of constant $\alpha=-1$ is that of a Killing vector field while the case of constant $\alpha=1$ is also of interest and applies for example to the conformal inflation of concentric spheres in $\R^3$.

\begin{proposition} Let $\mu,\nu\in C(M)$ obey the first order differential equations
\[ \tau(\mu)=\beta-(1+\alpha)\mu,\quad \tau(\nu)=\mu-\alpha\nu.\]
Then the calculus $\Omega^1(C(M)\rtimes_\tau\R)$ on  normal-ordered element $f(t)=\sum f_n t^n$ where $f_n\in C(M)$ (i.e. keeping the $t$-dependence to the right), obeys
\[ \theta' f(t)=f(t+\lambda\alpha)\theta',\quad \extd f=\bar\extd f+{\lambda \over 2}\theta'\bar\Delta f+ \del^0f\extd t+\lambda\Delta_0f\theta'\]
\[ \del^0 f(t)= {f(t)-f(t-\lambda)\over \lambda},\quad \Delta_0 f(t)= {\nu f(t+\lambda\alpha)+\mu  f(t-\lambda({\beta\over\mu}-\alpha))-(\nu+\mu)f(t+\lambda(\alpha-{\beta\over\nu+\mu}))\over \lambda^2}\]
and we also have
\[ [\extd t,f]=-\lambda\extd f+\lambda(\mu+\nu)({f(t+\lambda\alpha)-f(t+\lambda(\alpha-{\beta\over \mu+\nu}))\over\lambda})\theta'\]
\end{proposition}
\proof The behaviour on functions only on $M$ is already covered in Theorem~3.1. For a function purely of $t$ we prove the result at least for polynomials, by induction as follows (this generalises the bicrossproduct model case).  Assume $[\extd t,t^n]=p_n\extd t+q_n\theta'$. Then using the commutation relations,
\[p_n=(t-\lambda)p_{n-1}-\lambda t^{n-1},\quad q_n=(t-\lambda) q_{n-1}+\lambda\beta(t+\lambda\alpha)^{n-1}\]
which are solved by
\begin{equation}\label{dttn} [\extd t,t^n]=\left((t-\lambda)^n-t^n\right)\extd t +\mu\left((t+\lambda\alpha)^n-(t-\lambda({\beta\over\mu}-\alpha))^n\right)\theta'\end{equation}
provided $\mu$ obeys the relation stated. The proof for the $p_n$ is more elementary and left for the reader, while for $q_n$ we verify
that $q_1=\lambda\mu(\alpha + {\beta\over\mu}-\alpha)=\lambda\beta$ as required, and
\begin{eqnarray*} (t-\lambda) q_{n-1}&+&\lambda\beta(t+\lambda\alpha)^{n-1}\\
&=& \mu(t-\lambda(1+{\tau(\mu)\over\mu}))\left((t+\lambda\alpha)^{n-1}-(t-\lambda({\beta\over\mu}-\alpha))^{n-1}\right)+\lambda\beta(t+\lambda\alpha)^{n-1}=q_n
\end{eqnarray*}
taking account of the commutation relation $t \mu =\mu t-\lambda\tau(\mu)$. A further similar induction on $\extd t^n=t\extd t^{n-1}+[\extd t,t^{n-1}]+t^{n-1}\extd t$ provides the stated formulae as $\extd f=\del_0f+\lambda\Delta_0 f$.  Now suppose that $f=f(\ ,t)$ where the  dependence on $t$ is kept to the right and 
 combine the two cases via the Leibniz rule.  Note that with regard to the $t$-dependence $\theta'(\bar\Delta f)(t)=(\bar\Delta f)(t+\lambda\alpha)\theta'$ when our basic 1-forms are ordered to the right using the stated  commutation relation. Similarly we deduce from (\ref{dttn}) and the commutation $[\extd t,f]=-\lambda\extd f$ for a function on $M$ that in general for a normal ordered function
 \[ [\extd t, f(t)]=-\lambda\bd f(t)-{\lambda^2\over 2}\theta'\bDelta f(t)-\lambda\del_0 f(t)\extd t+\lambda \mu\left({f(t+\lambda\alpha)-f(t-\lambda({\beta\over\mu}-\alpha))\over\lambda}\right)\theta'\]
 which we can then write as stated. Note that there are also commutation relations between other differentials and functions. \endproof
 
 For example, if  $\tau(\alpha)=\tau(\beta)=0$ then $\mu=\beta/(1+\alpha)$ is killed by $\tau$ and  solves the $\mu$ equation. Similarly $\nu=\mu/\alpha$ is killed by $\tau$ and solves the $\nu$ equation. In this case $\nu+\mu=\beta/\alpha$ and 
\[\Delta_0 f=\beta{f(t+\alpha\lambda)+\alpha f(t-\lambda)-(1+\alpha)f(t)\over \lambda^2\alpha(1+\alpha)}.\]
If, moreover, $\alpha=1$ then we have $\Delta_0$ as $\beta\over 2$ times the standard symmetric finite difference Laplacian, while in the limit 
for $\alpha\to -1$ we have
\[\Delta_0 f{\buildrel \alpha\to -1\over{\  \longrightarrow\  }}{\beta\over\lambda}\left(\del_0 f- \dot f(t-\lambda)\right)\]
where $\dot f$ denotes the usual derivative in $t$, which is more readily seen to tend to ${\beta\over 2}\ddot f$ as $\lambda\to 0$. 

We see that the process of Theorem~3.1 induces a differential calculus in the extra `time' direction but it turns out to be of the finite-difference family that exists in noncommutative geometry even in one variable. Note also that $[t^n,\extd t]$ and $\extd t^n$  do not depend on the freedom in choices for $\mu,\nu$ as the boundary conditions and inductive relations do not depend on them, we only require them in order to have finite-difference type formulae, and we require them only locally, which is assured.  In particular, $\Delta_0$ does not depend on the choice of $\mu,\nu$ even if it looks as if it does, when the noncommutativity is taken into account.

\begin{corollary} Working in the calculus $\Omega^1(C(M)\rtimes\R)$, we define the {\em induced spacetime wave operator} $\square$ on $C(M)\rtimes_\tau\R$ by  $\extd f=\bd f+ \del_0f\, \extd t+{\lambda\over 2}\square f\, \theta'$.
\begin{enumerate}
\item    $\square f(t)=(\bDelta f)(t+\lambda\alpha)+2\Delta_0 f(t)$ on normal ordered $f(t)=\sum f_n t^n$.
\item In the classical limit  $ \lim_{\lambda\to 0} \square f= \bar\Delta f+ \beta\ddot f$ 
\item If  $\zeta^*=-{1\over 2}{\beta^{-1}}\bd \beta$ then the classical limit of $\square$ is the Laplace-Beltrami operator on $M\times\R$ for the static metric $\beta^{-1}\bd t\btens\bd t+\bar g$
\end{enumerate}
\end{corollary}
\proof  The shift by $\lambda\alpha$ in (1) is from the normal ordering. The classical limit (2) is a delicate computation assuming a Taylor expansion of $f$ about $t$ except that we have to be careful about the normal ordering. Writing $f=\sum_n f_n t^n$ we find
\begin{eqnarray*} \lambda^2\Delta_0f&=&\lambda^2 \sum_n f_n ({n\atop 2})(\nu(\alpha^2-\tau(\alpha)) +\mu((\alpha-{\beta\over \mu})^2-\tau(\alpha-{\beta\over\mu}))\\\
&&\quad\quad \quad\quad\quad-(\nu+\mu)((\alpha-{\beta\over\nu+\mu})^2-\tau(\alpha-{\beta\over\nu+\mu})))t^{n-2}+O(\lambda^3)\\
&=&\lambda^2{\beta\over 2}\ddot f+O(\lambda^3)\end{eqnarray*}
on computing $\tau$ as a derivation and using the defining equations for $\tau(\mu),\tau(\nu)$. The $\tau$ terms here arise from normal ordering of the different positions of $f$ in the linear term of the expansion
\[ (t+\lambda f)^n=t^n+ \lambda n f t^{n-1}-{n(n-1)\over 2}\lambda^2\tau(f)t^{n-2}+({n\atop 2})\lambda^2f^2 t^{n-2}+O(\lambda^3).\]
One can view this as the first terms of a noncommutative binomial identity for the action of a vector field on a function $f$. One can also derive the limit from the next lemma, but we have given the more direct proof. Part (3) is some elementary differential geometry. For a metric that splits as shown for some function $\beta\in C(M)$ in the time direction, the Levi-Civita connection is easily computed and $\bnabla\omega$ for the spacetime metric acquires an extra term from $\bd\beta$ and this enters into the spacetime time Laplace-Beltrami operator, which we compute as 
\[ \bar\square f=\beta\ddot f+\bDelta_{LB}-{1\over 2}\beta^{-1}\(\bd\beta,\bd f\).\]
 \endproof

This corollary makes good our philosophy that the `extra dimension' in the extended cotangent bundle expresses the Laplacian, as in Lemma~2.1, but now on the noncommutative spacetime version of $M\times\R$.  The following lemma provides more information about the time derivative component of the wave operator.

\begin{lemma} 
\[ \Delta_0 1=\Delta_0 t=0,\quad \Delta_0 t^2=\beta,\quad \Delta_0 t^3=3\beta t+\lambda((\alpha-1)\beta-2\tau(\beta))\]
while in general 
\[ \Delta_0 t^n=\sum_{i=0}^{n-1}(\del_0 t^{n-1-i})\beta(t+\lambda\alpha)^i,\quad [\extd t,t^n]=-(\del_0 t^n)\extd t+\lambda(\Delta_0 t^{n+1}-t\Delta_0 t^n)\theta'\]
\end{lemma}
\proof The first three cases are already contained in the proof of part 2 of Corollary~3.4 as there are no powers higher than $\lambda^2$ analysed there. For $t^3$ we note that
\begin{eqnarray*} (t+\lambda f)^3&=&t^3+\lambda( t^2 f+ t f t+ f t^2)+\lambda^2(f^2 t+f t f+ t f^2)+\lambda^3 f^3\\
&=&t^3+3\lambda f t^2+ 3\lambda^2 f^2 t+\lambda^3 f^3-3\lambda^2(\tau f)t+\lambda^3\tau^2 f-\lambda^3f\tau f-\lambda^3\tau( f^2)\end{eqnarray*}
We know in computing $\Delta_0 t^3$ from the expression in Proposition~3.3 that the order $1,\lambda$ terms wont contribute while the order $\lambda^3$ terms give us the classical contribution $3\beta t$. It remains to add up the terms at order $\lambda^3$, which contribute to $\Delta_0t^3$ the amount $\lambda\times$
\begin{eqnarray*} &&\nu(\alpha^3+\tau^2\alpha-3\alpha\tau\alpha)+\mu\left((\alpha-{\beta\over\mu})^3+\tau^2(\alpha-{\beta\over\mu})-3(\alpha-{\beta\over\mu})\tau(\alpha-{\beta\over\mu})\right)\\
&&\quad\quad-(\mu+\nu)\left((\alpha-{\beta\over\mu+\nu})^3+\tau^2(\alpha-{\beta\over\mu+\nu})-3(\alpha-{\beta\over\mu+\nu})\tau(\alpha-{\beta\over\mu+\nu})\right)\\
&=&3\alpha{\beta^2\over\mu}-3\alpha{\beta^2\over\mu+\nu}-{\beta^3\over\mu^2}+{\beta^3\over(\mu+\nu)^2}-\mu\tau^2({\beta\over\mu})+(\mu+\nu)\tau^2({\beta\over\mu+\nu})\\
&&\quad\quad+3(\alpha\mu-\beta)\tau({\beta\over\mu})-3(\alpha(\mu+\nu)-\beta)\tau({\beta\over\mu+\nu})
\end{eqnarray*}
which eventually simplifies to the result stated on repeated use of the relations
\[ \tau({\beta\over\mu})={\tau\beta\over\mu}-{\beta\over\mu}({\beta\over\mu}-(1+\alpha)),\quad  \tau({\beta\over\mu+\nu})={\tau\beta\over\mu+\nu}-{\beta\over\mu+\nu}({\beta\over\mu}-\alpha).\]
We omit the details in view of the general formula from which the final result can more easily be obtained. For the general formula, we compute
\[ \extd t^n=(\extd t^{n-1})t+t^{n-1}\extd t=(\del_0 t^{n-1})(\extd t)t+\lambda\Delta_0t^{n-1}\theta' t+t^{n-1}\extd t\]
\[
=(\del_0t^{n-1})t\extd t+\del_0 t^{n-1}\lambda(\beta\theta'-\extd t)+\lambda\Delta_0t^{n-1}(t+\lambda\alpha)\theta'+t^{n-1}\extd t\]
using the Leibniz rule and commutation relations. Comparing with $\extd t^n=(\del_0 t^n)\extd t+\lambda\Delta_0 t^n\theta'$ we deduce
\[ \del_0 t^n=(\del_0t^{n-1})(t-\lambda)+t^{n-1},\quad \Delta_0t^n=(\Delta_0t^{n-1})(t+\lambda\alpha)+\del_0t^{n-1}\beta.\]
The second of these provides the induction step easily solved to provide the result stated. We then use  the formula $\extd t^{n+1}=t\extd t^{n}+[\extd t,t^{n}]+t^{n}\extd t$ for the different induction in Proposition~3.3 now as a way to recover this $[\extd t,t^n]$ from $\Delta_0$.
 \endproof

We conclude with some immediate elements of the noncommutative geometry,  such as the natural extension to $\Omega^2$. We will not, however, take this too far in the present paper. Our main goal was the construction of the natural wave operator which we have done as part of $\Omega^1$ as covered in Proposition~3.3 and Corollary~3.4.

\begin{proposition}
\begin{enumerate}
\item  Natural relations in $\Omega^2(C(M)\rtimes_\tau\R)$ are provided by
\[ \{\theta',\extd t\}=-\lambda\extd \theta',\quad (\extd t)^2=-{\lambda\over 2}\extd(\beta\theta'),\quad [\extd\theta',t]=\lambda\extd\left((\alpha-1)\theta'\right)\]
\[ \{\omega,\extd t\}=-\lambda\extd\omega,\quad  \extd((\tau\beta+\beta(\alpha-1))\theta')=0\]
for all $\omega\in\bar\Omega^1$, in addition to those in Proposition~\ref{omega2}, (\ref{thetawedge}), (\ref{thetaf}).
\item Given $\mu,\nu$ as in Proposition~3.3, the calculus is inner if the last condition in (1) is replaced by the stronger $\extd((\beta-(\mu+\nu))\theta')=0$. Here 
\[ \theta:=\extd t-(\mu+\nu)\theta',\quad [\theta,\ \}=-\lambda\extd \]
where $[\ ,\ \}$ denotes (anti)commutator on degree 0 (degree 1) respectively. \end{enumerate}
\end{proposition}
\proof  We start with the maximal prolongation: we apply $\extd$ to the relations in degree 1 involving $t$ to obtain
\begin{equation}\label{maxrel1} \{\extd f,\extd t\}=0,\quad (\extd t)^2=-{\lambda\over 2}\extd(\beta\theta'),\quad \{\theta',\extd t\}=[\extd\theta',t]-\lambda\extd(\alpha\theta')\end{equation}
for all $f\in C(M)$. Using Lemma~3.5 we find no further restriction from $\extd^2t^2=0$. Thus, 
\[ \extd^2t^2=\extd((2t-\lambda)\extd t+\lambda\beta\theta')=2(\extd t)^2+\lambda\extd(\beta\theta')=0\]
is one of the relations already found. Next from Lemma~3.5,
\begin{eqnarray*}
\extd^2 t^3&=&\extd((3t^2-3\lambda t+\lambda^2)\extd t+\lambda(3t\beta+\lambda\tau\beta+\lambda\beta(\alpha-1))\theta')\\
&=&(6(t-\lambda)(-{\lambda\over 2}\extd(\beta\theta'))+3\lambda\beta\theta'\extd t+3\lambda(\extd t)\beta\theta'+\lambda 3t\extd(\beta\theta')+\lambda^2\extd((\tau\beta+\beta(\alpha-1))\theta')\\
&=&\lambda^2\extd((\tau\beta+(\alpha-1)\beta)\theta')+3\lambda\beta(\{\theta',\extd t\}+\lambda\extd\theta').
\end{eqnarray*}
Hence we require
\begin{equation}\label{maxrel2} 3\beta(\{\theta',\extd t\}+\lambda\extd\theta')+\lambda\extd((\tau\beta+(\alpha-1)\beta)\theta')=0
\end{equation}
for $\extd t^3=0$. For  $\extd t^4$ we have (using Lemma~3.5 for the second)
\[ \del_0 t^4=(4 t^3-6\lambda t^2+4\lambda^2 t-\lambda^3)\extd t\]
\begin{eqnarray*} \Delta_0t^4&=&\left(6 t^2\beta+4\lambda t(\tau\beta+\beta(\alpha-1))+\lambda^2(-\tau\beta-\beta(\alpha-1)+\tau^2\beta+\tau(\alpha\beta)+(\tau\beta)\alpha+\beta\alpha^2)\right)\theta'\\
&=&\left(6 t^2\beta+4\lambda t \, h_2+\lambda^2(\tau h_2+\alpha h_2-\beta)\right)\theta'\\
\end{eqnarray*}
where $h_2=\tau\beta+(\alpha-1)\beta$, say. From this and a lengthy computation using (\ref{maxrel1})-(\ref{maxrel2}) we find
similarly
\[\extd^2t^4=\extd(\del_0 t^4\extd t+\lambda\Delta_0 t^4\, \theta')=\lambda^2(4h_2+3\beta)(\{\theta',\extd t\}+\lambda\extd\theta')+\lambda^3\extd((\tau h_2+(\alpha-1)h_2)\theta')=0 \]
as one of a sequence of identities obtained by repeatedly applying $[\ ,t]$ to (\ref{maxrel2}). To see this let $D_\tau(f)=\tau f+(\alpha-1)f$ for any function $f$ and $h_1=\beta$, $h_i=D_\tau(h_{i-1})$. Then expanding out and repeatedly applying the commutation relations we find  that
\[ [f(\{\theta',\extd t\}+\lambda\extd\theta'),t]=\lambda(f+D_\tau f)(\{\theta',\extd t\}+\lambda\extd \theta'),\quad
[\extd(f\theta'),t]=f(\{\theta',\extd t\}+\lambda\extd\theta')+\lambda\extd ((D_\tau f)\theta')\]
for any function $f$. Setting $f=\beta$ gives the identity needed. Although we have not written out a formal proof for general $n$ it is clear that there are no further relations from $\extd^2 t^n=0$ beyond (\ref{maxrel2}). 

So far we have only prolonged from the commutation relations involving $t$ in (\ref{cr}); if we want to have a tensorial construction coming from the manifold $M$ it is natural to restrict further. In particular, from the first of (\ref{maxrel1}) we find
\begin{eqnarray*}\{a_i\bd b_i,\extd t\}&=&a_i\{\bd b_i,\extd t\}-[a_i,\extd t]\bd b_i=-{\lambda\over 2}a_i\{\bDelta b_i\theta',\extd t\}
-\lambda\extd a_i\, \bd b_i\\
&=&-{\lambda\over 2}a_i\bDelta b_i\{\theta',\extd t\}+{\lambda^2\over 2}a_i\bd \bDelta b_i\theta' 
-\lambda\extd a_i\, \bd b_i\\
&=&-\lambda\extd(a_i\bd b_i)-{\lambda\over 2}a_i\bDelta b_i(\{\theta',\extd t\}+\lambda\extd\theta')
\end{eqnarray*}
using $\extd^2 b_i=0$. It is then natural to impose the first stated relation in part (1)  in order that the right hand side depends only on $\omega=a_i\bd b_i \in \bar\Omega^1$. This then yields the remaining relations.

For (2), if we seek $\theta=\extd t-g\theta'$ with the required `inner' property for some function $g$ then we will need $\lambda\extd t=[t,\extd t-g\theta']=\lambda\extd t+\lambda\beta\theta'+\tau(g)\theta'+g\lambda\alpha\theta'$, which requires $g$ to obey the equation $\tau g+\alpha g=\beta$, which in turn is solved by $g=\mu+\nu$ in Proposition~3.3. We also require $-\lambda\extd \theta'=\{\theta',\extd t-g\theta'\}=\{\theta',\extd t\}$ independently of $g$, so this is a natural restriction also from the point of view of having an inner calculus with a form of $\theta$ close to $\extd t$. Finally, we require $0=\{\extd t,\theta\}=\{\extd t,\extd t-g\theta'\}=-\lambda\extd(\beta\theta')-g(\{\theta',\extd t\}+\lambda\extd\theta')+\lambda\extd (g\theta')$. In our context this leads to the relation stated in (2). Note, in general, however, that, 
\[ 0=g(\{\theta',\extd t\}+\lambda\extd\theta')+\lambda\extd((\beta-g)\theta')\]
implies by application of $[\ ,t]$ and the identities in part (1) that 
\[ 2\beta (\{\theta',\extd t\}+\lambda\extd\theta')+\lambda\extd( h_2\theta')=0.\]
Comparison with (\ref{maxrel2}) again makes it natural to restrict to $\{\theta',\extd t\}+\lambda\extd\theta'=0$. In this case the last of the relations stated in (1) is implied by our simpler requirement $\extd((\beta-g)\theta')=0$ for the calculus to be inner in the form taken. Finally, we already have  $[f ,\extd t]=\lambda\extd f$ on $f\in C(M)$ and  $\{\omega,\extd t\}=-\lambda\extd\omega$ (obtained in part (1)) and these are not affected by the addition of any functional multiple of $\theta'$ as $\theta'^2=\{\omega,\theta'\}=0$.  \endproof

Note that locally we can always solve for $\mu+\nu$ as required for an inner calculus while taking  $\extd((\beta-(\mu+\nu))\theta')=0$ typically as a definition of $\extd\theta'$. In that sense one can say that the calculus on $C(M)\rtimes_\tau\R$ is always `locally inner'   at least up to degree 2. When $\beta,\alpha\ne 0$ are constant, for example, we have $\beta-(\mu+\nu)=\beta(\alpha-1)/\alpha$ so we need $\extd\theta'=0$ unless $\beta=0$ or  $\alpha=1$. The property of being inner is a desirable feature of any sufficiently noncommutative calculus and in this sense we see that the calculus is better behaved than in Section~2 before we adjoined $t$. One can use the requirements for $\theta$ to similarly define $\extd$ and find relations in all degrees, but we defer this to a sequel. We conclude with a further miscellaneous observation.

\begin{proposition}
For any $h\in C(M)$ the transformation 
\[ \extd t\to \extd' t=\extd t+h\theta',\quad \beta\to \beta'=\beta+\tau h+(\alpha+1)h\]
 gives a differential calculus $(\Omega^1,\extd')$ with isomorphic bimodule structure to that of $(\Omega^1,\extd)$.
 \end{proposition}
 \proof Let $(\Omega,\extd')$ be the calculus constructed as in Theorem~3.1 but with parameter $\beta'$ in place of $\beta$. The map $\phi:(\Omega^1,\extd)\to( \Omega^1,\extd')$ defined by $\phi(\extd t)=\extd' t-h\theta'$ is a bimodule map. To see this, we check consistency with the bimodule commutation relations involving $\extd t$. Thus  $\phi([\extd t,t])=[\extd' t-h\theta',t]=\lambda\beta'\theta'-\lambda\extd' t-h\lambda\alpha\theta'-\lambda\tau(h)\theta'=\lambda\beta\theta'-\lambda(\extd' t-h\theta')=\phi(\lambda\beta\theta'-\lambda\extd t)$. Similarly $\phi([f,\extd t])=[f,\extd' t-h\theta']=\phi(\lambda\extd t)$.   It is important to note, however, that $\phi$ does not form a commutative triangle connecting $\extd,\extd'$, i.e. these are {\em not} necessarily isomorphic as differential calculi.   From Lemma~3.5 we see rather that
 \[ \phi(\extd t^n)=\extd' t^n-\left(\del_0t^nh+\lambda\sum_{i=0}^{n-1}\del_0t^{n-1-i}(\tau h+(\alpha+1)h)(t+\lambda\alpha)^i\right)\theta'\]
  showing how the differentials change when there is $t$-dependence. Meanwhile, $\phi(\extd f)=\extd' f$ for all $f\in C(M)$.
 \endproof

The first order nature of the differential on $h$ here means that, at least locally, any two choices of $\beta$ have the same bimodule structure up to isomorphism, i.e this aspect of the structure in Theorem~3.1 is not being changed by the variation of $\beta$, up to isomorphism. For example $h=-\mu$ renders $\beta'=0$. Globally, however, we expect nontrivial equivalence classes of calculi, in the sense of isomorphism classes of bimodule structures, depending on the topology of $M$. 

\section{Quantisation of the Schwarzschild black hole and other spherically symmetric static spacetimes}

In this section we apply the machinery above to construct, in principle, the quantum black hole differential algebra.  Most of the work is in the preliminary Section 4.1 is to give an entirely algebraic approach to polar coordinates in classical Riemannian geometry. We are then in position in Seciton 4.2 to read off the quantum version.

\subsection{Algebraic polar coordinates and the Schwarzschld solution}

Most of the section covers familiar ground, except that we do so as a novel application of algebraic methods developed for noncommutative geometry. Unless stated otherwise, we work over a general field $k$ of characteristic not 2. And because everything in this subsection is classical, we will endeavour to put bars over the geometrical objects. Where omitted for brevity, they should be understood.

We start by considering algebraic analogues of the space $\R^3\setminus\{0\}$ where we delete the origin. We assume that we work in some sufficiently large `coordinate algebra' $A$ containing at least  the mutually commuting generators $x_i,r$ for the Cartesian and radial coordinates subject to the relation $r^2=\sum_i x_i^2$, as well as $r^{-1}$ and sufficiently many other rational functions of $r$ as to be able to solve any equations we need.  We can think of  the coordinate algebra $A=k[x_1,x_2,x_3][(r)$ (rational functions in $r$ and polynomials in the $x_i$) subject to the above relation. We consider
\begin{equation}\label{proj}e_{ij}=\delta_{ij}-{x_ix_j\over r^2},\quad e\in M_3(A)\end{equation}
which is a projector,
\[ e^2_{ik}=e_{ij}e_{jk}=\delta_{ik}+{x_ix^2x_k\over r^4}-2{x_ix_k\over r^2}=e_{ik};\quad x_ie_{ij}=0.\]
 According to standard arguments in (commutative and) noncommutative geometry 
 $\CE=eA^{\oplus 3}$ is a projective module or `vector bundle' and as such it has on it a canonical nontrivial `Grassmann connection' $\bnabla$ at least for the universal calculus. Thus far we do not actually need the generators to commute. However, we now assume that they do and that there is a standard commutative differential structure, thus $\{\bd x_i\}$ define a free module of 1-forms $\bar\Omega^1$  with    $\bd r^2=2r\bd r=2x_i\bd x_i$. Then $\CE$ is a 2-dimensional sub-bundle spanned by
 \[ \omega_i=e_{ik}\bd x_k= \bd x_i-{x_i \bd r\over r},\quad x_i\omega_i=0\]
as tangential to the sphere about the origin passing through any point. We have  $\bar\Omega^1(A)=\CE\oplus A\bd r$ as a direct sum of tangential and radial bundles. The Euclidean metric can be written as
\[\bar\eta=\bd x_i \btens\bd x_i=(\bd x_i-{x_i \bd r\over r})\btens(\bd x_i-{x_i \bd r\over r})+\bd r\btens\bd r=\omega_i\btens\omega_i+\bd r\btens\bd r\]
so $\bar\eta^{ang}=\omega_i\btens\omega_i$ provides the `angular part' (usually denoted $r^2\bd\bar\Omega$) of the metric in polar coordinates but written in tensor form.

\begin{proposition}
The Grassmann connection for the projector $e$ is
\[  \bnabla\omega_i=-{x_i\over r^2}\omega_j\btens\omega_j ,\]
 and is metric compatible with $\bar\eta^{ang}$ but not torsion free and has curvature
\[ R_\bnabla(\omega_i)=-{ \omega_i\over r^2}\omega_j\btens \omega_j.\]
It extends to a connection on $\bar\Omega^1$ by $\bnabla(\bd r)=0$ which is metric compatible with $\bar\eta$ but not torsion free and has $R_\bnabla(\bd r)=0$.
\end{proposition}
\proof The connection on a projective module is provided by in our application  by $\bnabla\omega=(\bd e)e \btens\omega$. We compute
\[ (\bd e_{ij}) e_{jk}=-\bd({x_ix_j\over r^2})e_{jk}=-{x_i\over r^2}\bd x_j e_{jk}=-{x_i\over r^2}\omega_j e_{jk}.\]
It follows that $\bnabla$ is well-defined but one can also verify it directly
\[ \bnabla(x_i\omega_i)=-x_i\omega_j\btens\omega_j{x_i\over r^2}+\bd x_i\btens\omega_i=-\omega_j\btens\omega_j+\omega_i\btens\omega_i=0.\]
The curvature can be computed similarly from $-(\bd e\bd e)e$ but we prefer to compute it directly 
\begin{eqnarray*} R_{\bnabla}(\omega_i)&=&(\bd\btens\id-(\wedge\btens\id)(\id\btens\bnabla))\bnabla\omega_i=(\bd\btens\id-(\wedge\btens\id)(\id\btens\bnabla))(-{ x_i\over r^2}\omega_j \btens\omega_j)\\
& =&-\bd({ x_i\over r^2})\omega_j\btens\omega_j-{ x_i\over r^2}\bd(\omega_j)\btens\omega_j+{ x_i\over r^2}\omega_j\wedge\bnabla\omega_j\\
&=&-({\omega_i\over r^2}-{x_i\over r^3}\bd r)\omega_j\btens\omega_j+{ x_i\over r^2}\omega_j {\bd r\over r}\btens\omega_j+0
\end{eqnarray*}
since $\omega_ix_i=0$. Assuming antisymmtry in the exterior product gives the result as stated.
Similarly, the torsion tensor $T_\bnabla$ and metric compatibility are
\[ T_\bnabla(\omega_i)=\bnabla\wedge\omega_i-\bd \omega_i=-\bd\omega_i=\bd\left({x_i\over r}\right)\bd r=(\bd x_i){\bd r\over r}-{x_i \over r^2}\bd r\bd r=\omega_i{\bd r\over r}\]
\[\bnabla(\bar\eta^{ang})= \bnabla(\omega_i\btens\omega_i)=-\omega_j\btens\omega_j {x_i\over r^2}\btens \omega_i-\omega_j \btens\omega_i\btens \omega_j {x_i\over r^2}=0\]
where the left output of the action of $\bnabla$ is kept to the far left. We have assumed familiar properties of classical connections. \endproof 

This is a type of monopole on $\R^3\setminus\{0\}$ but can also be viewed as a Riemannian connection with torsion as explained. It restricts on spheres of constant radius (where $\bd r$ is projected out) to the Levi-Civita connection on these. 
We compare this connection with the more obvious trivial connection. This is defined by $\bnabla(\bd x_i)=0$ and is obviously torsion free and metric compatible with $\bar\eta$ and has $R_{\bnabla}=0$.

\begin{proposition} The trivial connection  can be written in polar coordinates as 
\[ \bnabla(\omega_i)=-{x_i\over r^2}\omega_j\btens\omega_j-{\omega_i\over r^2}\btens r\bd r,\quad \bnabla(\bd r)={1\over r}\omega_i\btens\omega_i\]
\end{proposition}
\proof For the trivial connection $\bnabla(\bd x_i)=0$ we compute
\[ \bnabla(r\bd r)=\bd r\btens\bd r+r\bnabla(\bd r)=\bnabla(x_i\bd x_i)=\bd x_i\btens\bd x_i+0=\bar\eta\]
from which we conclude  $\bnabla(\bd r)$ as stated. Also
\[ \bnabla(\omega_i)=0-\bnabla({x_ir\over r^2}\bd r)=-\bd({x_i\over r^2})\btens r\bd r-{x_i\over r^2}\bnabla(r\bd r) \]
\[{\ }\quad\quad =-{\bd x_i\over r^2}\btens r\bd r+{2 x_i\over r^4}r\bd r\btens r\bd r-{x_i\over r^2}\bd x_j\btens\bd x_j\]
which provides $\bnabla(\omega_i)$ as stated. One can verify directly that these expressions give a (flat) torsion free and metric compatible connection as they must. Indeed, the second term of $\bnabla(\omega_i)$ precisely kills the torsion compared to Proposition~4.1 but now introduces
\[ \bnabla(\bar\eta^{ang})=-{\omega_i\over r^2}\btens r\bd r\btens\omega_i-{\omega_i\over r^2}\btens\omega_i\btens r\bd r\]
compared to the calculation in Proposition~4.1. This is precisely compensated by $\bnabla(\bd r\btens\bd r)$ as now $\bnabla(\bd r)\ne0$, allowing the connection to remain metric compatible for $\bar\eta$. Finally, we still have $\omega_j\bnabla(\omega_j)=0$ so the previous calculation for $R_\bnabla(\omega_i)$ just acquires extra terms from the new part of $\bnabla(\omega_i)$:
\[R_\bnabla(\omega_i)=-{ \omega_i\over r^2}\omega_j\btens \omega_j -\bd({\omega_i\over r^2})\btens r\bd r+{\omega_i\over r^2} \bnabla(r\bd r)=-{ \omega_i\over r^2}\omega_j\btens \omega_j-{\omega_i\over r^2}\bd r\btens\bd r+{\omega_i\over r^2} \bnabla(r\bd r)\]
which indeed cancel to give 0 from the form of $\bnabla(r\bd r)$. In addition
\[ R_\bnabla(\bd r)=(\bd\btens\id-(\wedge\btens\id)(\id\btens\bnabla))({1\over r}\omega_i\btens\omega_i)=\bd({\omega_i\over r})-{\omega_i\over r}\bnabla(\omega_i)=0\]
as both terms vanish. Of course we already know the results here but these computations will be a model for the next proposition. \endproof

Now looking carefully at the mentioned direct check of how the trivial $\bnabla$ gets to be torsion free and metric compatible in the above `radial/tangential' framework,  one can see that the proof can be generalised to the following:

\begin{proposition} Let 
\[ \bar g=h(r)^2\bd r\btens\bd r+ \omega_i\btens \omega_i \]
for a function $h(r)$, which we assume to be invertible in the algebra. Then
\[ \bnabla(\omega_i)=-{x_i\over r^2}\omega_j\btens\omega_j-{\omega_i\over r^2}\btens r\bd r,\quad \bnabla(\bd r)={1\over h(r)^2 r}\omega_i\btens\omega_i-{h'(r)\over h(r)}\bd r\btens\bd r\]
is torsion free, metric compatible with $\bar g$ and has curvature
\[ R_\bnabla(\omega_i)=-{\omega_i\over r^2}\left(\left(1-{1\over h(r)^2}\right)\omega_j\btens\omega_j+{h'(r) r\over h(r)}\bd r\btens\bd r\right),\quad R_\bnabla(\bd r)=-{h'(r)\over h(r)^3 r}\bd r\, \omega_i\btens\omega_i.\]
\end{proposition}
\proof Because the connection on `angular' forms is unchanged, this part remains torsion free. For the radial part clearly $\bnabla\wedge\bd r=0$ so the torsion on $\bd r$ also vanishes. We have to check metric compatibility and we note 
that
 \[ \bnabla(h(r) \bd r)={1\over h(r) r}\omega_i\btens\omega_i\]
similarly to the structure in Proposition~4.2. The computation of $\bnabla(\omega_i\btens\omega_i)$ is unchanged and now clearly killed in just the same way by $\bnabla(h(r)\bd r\btens h(r) \bd r)$. It remains to compute the curvature. The only difference for $R_\bnabla(\omega_i)$ compared to the direct calculation in Proposition~4.2 is the form of $\bnabla(r\bd r)$, so this time
\[ R_\bnabla(\omega_i)=-{\omega_i\over r^2}\bar\eta+{\omega_i\over r^2}\left({1\over h(r)^2}\omega_j\btens\omega_j+(1-{h'(r)\over h(r)}r)\bd r\btens\bd r\right)\]
giving the result as stated. Similarly
\[ R_\bnabla(h(r)\bd r)=(\bd\btens\id-(\wedge\btens\id)(\id\btens\bnabla))({1\over h(r) r}\omega_i\btens\omega_i)=\bd({1\over h(r)}){\omega_i\over r}\btens\omega_i\]
giving the stated result.
\endproof

Incidentally, the connection on radial forms can also be written more compactly as
 \begin{equation}\label{nablahr} \bnabla(h(r) r\bd r)={1\over h(r)}\bar g.\end{equation}

\begin{corollary} The standard lifting $i:\bar\Omega^2\to \bar\Omega^1\btens\bar\Omega^1$ and trace applied to the curvature in Proposition~4.3 gives 
\[ {\rm Ricci}=-{1\over 2r}\left({h'(r)\over  h(r)^3 }+{1\over r}\left(1-{1\over h(r)^2}\right) \right)\omega_j\btens\omega_j-{h'(r) \over h(r)r}\bd r\btens\bd r.\]
In particular,  ${\rm Ricci}\propto \bar g$ (an Einstein space) iff $r h'(r)= h(r)(h(r)^2-1)$. \end{corollary}
\proof Strictly speaking we repeat the computations for ${\rm Ricci}$ from the curvature in Proposition~4.3, but the result is the same as setting all terms involving $f$ to zero in the preceding theorem and hence we omit the details. Hence to be an Einstein space we need 
\[ {h'(r)\over  h(r)^3 }+{1\over r}\left(1-{1\over h(r)^2}\right)=2 {h'(r)\over  h(r)^3 }\]
as stated. \endproof

Over $\R$ the equation is solved by $h(r)=1/\sqrt{1+K r^2}$ where $K$ is a parameter. Then ${\rm Ricci}=K \bar g$. Hence this is space of constant curvature. For $K<0$ (in our conventions) it is essentially $S^3$, while for $K>0$ it is hyperbolic 3-space. In both cases the removal of $r=0$ is not required other than for use of our polar coordinates. Note that we cannot by contrast solve ${\rm Ricci}=0$ unless $h(r)=1$ which is the case of Proposition~4.2. 

Next, we consider Killing vector fields needed later. Because in the present section we prefer differential forms we map a vector field $\tau$ to a 1-form $\tau^*=\bar g(\tau)$. In these terms a conformal Killing 1-form is required to obey 
\[ (\id+\sigma)\bnabla\tau^*\propto \bar g\]
where $\sigma$ is the trivial flip map.  

\begin{corollary} The metric in Proposition~4.4 has conformal Killing forms $\tau^*$ and Killing forms $\tau_i^*$  (not linearly independent),
\[ \tau^*=h(r)r\bd r,\quad \tau^*_i=x_j\omega_k \eps_{ijk}\]
where $\eps_{ijk}$ is the totally antisymmetric tensor with $\eps_{123}=1$. 
\end{corollary}
\proof That $\tau^*$ is a conformal Killing form is immediate from (\ref{nablahr}), which indeed says that $\bnabla\tau^*=\bar g/h(r)$. If one tries a more general form $\tau^*=f(r)h(r)r\bd r$ then one can deduce that $f'(r)=0$. The $\tau^*_i$ correspond to the action of the group of rotations in 3-dimensions and one easily computes 
\[ \bnabla\tau_i^*=\eps_{ijk}(\omega_j\btens\omega_k+{x_j\over r}(\bd r\btens\omega_k-\omega_k\btens\bd r))\]
which is manifestly antisymmetric. Hence its symmetrization vanishes and we have a Killing 1-form.  \endproof

We mention that at least over $\R$ all of the computations in this section work equally well with vector fields and the inverse metric. Some natural vector fields for our `polar coordinates' are
\[ \rho = x_i\del^i=r{\del\over\del r},\quad e_i=\del^i - {x_i\over r^2}\rho;\quad x_ie_i=0\]
where $\del^i={\del\over\del x_i}$ are the vector fields for the Cartesian coordinates. Note that $\rho$ acts as the degree operator so $\rho(x_i)=x_i$ and that the $e_i$ are not linearly independent. The latter are partial derivatives associated to the $\omega_i$ in the sense
\[ \bd \psi=({\del\over\del r}\psi)\bd r+e_i(\psi)\omega_i\]
for all functions $\psi$ in our coordinate algebra $A$ on $\R^3\setminus\{0\}$. Indeed, we have by easy computations
\[ \<\rho,\bd r\>=r,\quad \<\rho,\omega_i\>=0,\quad \<e_i,\bd r\>=0,\quad \<e_i,\omega_j\>=e_{ij} \]
where we use our projector matrix entries and in this sense also the $e_i$ are `dual' to the $\omega_i$. Next, we view $\bar g$ as a map  $\bar\Omega^{-1}\to \bar\Omega^1$ by evaluation against the first factor (say) and as an application, we compute for the inverse of the metric in Proposition~4.3,
\[ \bar g(\rho)=h(r)^2 r\bd r,\quad \bar g(e_i)=e_{ij}\omega_j=\omega_i; \quad \bar g^{-1}(\bd r)={1\over h(r)^2 r}\rho,\quad \bar g^{-1}(\omega_i)=e_i\]
Hence
\begin{equation}\label{invmetric} \(\bd r,\bd r\)={1\over h(r)^2},\quad \(\bd r,\omega_i\)=0,\quad \(\omega_i,\omega_j\)=e_{ij}.\end{equation}
Then, for example, the corresponding conformal Killing vector fields in Corollary~4.5 are 
\begin{equation}\label{taui} \tau={1\over h(r)}\rho,\quad \tau_i={\eps_{ijk}\over h(r)^2}x_je_k.\end{equation}

Finally, we can go one step further with the above as the spatial part of a radially-symmetric static spacetime geometry. Thus, we consider $A[t]$ with a new variable (time) adjoined as the coordinate algebra of $(\R^3\setminus\{0\})\times\R$.  We assume that $\bar\Omega^1(A[t])$ is spanned as before but with the additional 1-form $\bd t$, and that everything remains commutative. 

\begin{proposition}  Let 
\[ \bar g=-f(r)^2\bd t\btens\bd t+h(r)^2\bd r\btens\bd r+ \omega_i\btens \omega_i\]
for functions $h(r),f(r)$, which we assume to be invertible in the algebra. Then $\bnabla(\omega_i)$ as before and 
\[  \bnabla(\bd r)={1\over h(r)^2 r}\omega_i\btens\omega_i-{h'(r)\over h(r)}\bd r\btens\bd r-{f'(r)f(r)\over h(r)^2}\bd t\btens\bd t,\quad \bnabla(\bd t)=-{f'(r)\over f(r)}(\bd t\btens\bd r+\bd r\btens\bd t)\]
is torsion free, metric compatible with $\bar g$ and has curvature
\[ R_\bnabla(\omega_i)=-{\omega_i\over r^2}\left(\left(1-{1\over h(r)^2}\right)\omega_j\btens\omega_j+{h'(r) r\over h(r)}\bd r\btens\bd r+{f'(r)f(r)r\over h(r)^2}\bd t\btens \bd t \right)\]
\[ R_\bnabla(\bd r)=-{h'(r)\over h(r)^3 r}\bd r\, \omega_i\btens\omega_i+{f(r)\over h(r)^3}\left(f'(r)h'(r)-f''(r)h(r)\right)\bd r\bd t\btens\bd t\]
\[ R_\bnabla(\bd t)={1\over f(r)h(r)}\left(f'(r)h'(r)-f''(r)h(r)\right)\bd r\bd t\btens\bd r+ {f'(r)\over f(r)h(r)^2 r}\bd t\, \omega_i\btens\omega_i.\]
\end{proposition}
\proof The torsion on $\bd r$ continues to vanish as $(\bd t)^2=0$ and vanishes on $\bd t$ by $\{\bd t,\bd r\}=0$; so the connection remains torsion free. For metric compatibility we write the connection in the form
 \[ \bnabla(h(r) \bd r)={1\over h(r) r}\omega_i\btens\omega_i-{f'(r)f(r)\over h(r)}\bd t\btens\bd t,\quad \bnabla(f(r)\bd t)=-f'(r)\bd t\btens\bd r\]
 Then compared to the previous case $\bnabla(h(r)\bd r\btens h(r)\bd r)$ acquires an extra term
 \[ -{f'(r)f(r)\over h(r)}(\bd t\btens\bd t\btens h(r)\bd r+\bd t \btens h(r)\bd r\btens \bd t)\]
 which is exactly cancelled by $\bnabla(-f(r)\bd t\btens f(r)\bd t)$. Hence the connection remains metric compatible. As $\bnabla(\omega_i)$ is unchanged, $R_\bnabla(\omega_i)$ in the previous computation is affected only through $\bnabla(r\bd r)$ which acquires an extra $-(f'(r)f(r)/h(r)^2)r\bd t\btens\bd t$, giving the additional contribution stated. Similarly, in the previous computation of $R_\bnabla(h(r)\bd r)$ the change in $\bnabla(h(r)\bd r)$ gives an additional contribution
 \[(\bd\btens\id-(\wedge\btens\id)(\id\btens\bnabla))(-{f'(r)f(r)\over h(r)}\bd t\btens\bd t)=-\bd \left({f'(r)f(r)\over h(r)}\right)\bd t\btens\bd t+{f'(r)f(r)\over h(r)}\bd t\wedge\bnabla\bd t\]
 giving the additional contribution stated. Finally, we compute
 \[ R_\bnabla(f(r)\bd t)=(\bd\btens\id-(\wedge\btens\id)(\id\btens\bnabla))(-f'(r)\bd t\btens \bd r)\]
 \[{\ }\quad\quad\quad\quad=-\bd(f'(r))\bd t\btens\bd r+f'(r)\bd t\wedge\left({1\over h^2 r}\omega_i\btens\omega_i-{h'(r)\over h(r)}\bd r\btens\bd r\right)\]
 as $(\bd t)^2=0$. This gives the result stated. 
  \endproof

\begin{theorem} The standard lifting $i:\bar\Omega^2\to \bar\Omega^1\btens\bar\Omega^1$ and trace applied to the curvature in Proposition~4.3 gives 
\begin{eqnarray*} {\rm Ricci}&=&{1\over 2r}\left( {f'(r)\over  f(r)h(r)^2 }-{h'(r)\over  h(r)^3 }-{1\over r}\left(1-{1\over h(r)^2}\right) \right)\omega_j\btens\omega_j\\
&&-\left({1\over 2f(r)h(r)}\left(f'(r)h'(r)-f''(r)h(r)\right)+{h'(r) \over h(r)r}\right)\bd r\btens\bd r\\
&&+\left({f(r)\over 2 h(r)^3}\left(f'(r)h'(r)-f''(r)h(r)\right)-{f'(r)f(r)\over h(r)^2r}\right)\bd t\btens\bd t. \end{eqnarray*}
In particular, ${\rm Ricci}=0$ if $h(r)=1/f(r)$ and $r{d\over d r}f(r)^2=1-f(r)^2$. \end{theorem}
\proof The standard lift of 2-forms in classical geometry is to identify them with antisymmetric tensors, so for example $i(\bd r\bd t)={1\over 2}(\bd r\btens\bd t-\bd t\btens\bd r)$. We then take a trace of $(i\btens\id)R_\bnabla$ as an operator mapping to the first tensor factor (say) of its output, to give ${\rm Ricci}$. When doing this, clearly $R_\bnabla(\bd r)\propto \bd r \bd t\btens\bd t$ will contribute ${1\over 2}\bd t\btens\bd t$ to the trace as only the first term of the lift will contribute.  Similarly for the contribution from $R_\bnabla(\bd t)$. For $R_\bnabla(\omega_i)$, where the a term is of the form $\omega_i X$ and $X$ does not involve $\{\omega_j\}$ in its first tensor factor, we will similarly have $(X/2)\times 2=X$ for the contribution to the trace from $\omega_i\to {1\over 2}\omega_i\btens X$, because  the projective module has rank $2$ and the operation is as a multiple of the identity. For a term in $R_\bnabla(\omega_i)$ of the form $\omega_i\omega_j\btens\omega_j$ will again have $\omega_j\btens\omega_j$ for the same reason but also the trace of 
\[ \omega_i\mapsto -{1\over 2}\omega_j\btens\omega_i\btens\omega_j\]
from the antisymmetrisation. This will contribute $-{1\over 2}\omega_j\btens\omega_j$ giving a total contribution from such a term in $R_\bnabla(\omega_i)$ of ${1\over 2}\omega_j\btens\omega_j$. With these observations, we see without further computation that 
\begin{eqnarray*} {\rm Ricci}&=&-{1\over r^2}\left({1\over 2}\left(1-{1\over h(r)^2}\right)\omega_j\btens\omega_j+{h'(r) r\over h(r)}\bd r\btens\bd r+{f'(r)f(r)r\over h(r)^2}\bd t\btens \bd t \right)\\
&&-{h'(r)\over 2 h(r)^3 r} \omega_i\btens\omega_i+{f(r)\over 2 h(r)^3}\left(f'(r)h'(r)-f''(r)h(r)\right)\bd t\btens\bd t\\
&&-{1\over 2f(r)h(r)}\left(f'(r)h'(r)-f''(r)h(r)\right)\bd r\btens\bd r+ {f'(r)\over 2 f(r)h(r)^2 r} \omega_i\btens\omega_i.\end{eqnarray*}
which then combines as stated. Note that if ${\rm Ricci}=0$ then combining the $\bd t\btens\bd t$ and $\bd r\btens\bd r$ equations we deduce that $f'/f+h'/h=(fh)'/(fh)=0$ which over $\R$ implies that $h\propto 1/f$, and so on, but this depends on the field.
 \endproof

Over $\R$ the equation is solved by $f(r)=\sqrt{1-{\gamma\over r}}$ where $\gamma$ is a parameter. This is the Schwarzschild black hole with event horizon at $r=\gamma$. Note that our natural algebraic conventions in defining ${\rm Ricci}$ differ by $\pm{1\over 2}$ from the more usual ones used in other sections of the paper but this does not affect Ricci flatness of course.  Also note that it was convenient (and conventional in physics) but  not essential to work with $f,h$ -- one can work just as well throughout with $f^2(r)$ and $h^2(r)$ as the functions of interest. All formulae can be reworked in terms of these without square roots and one should do so for a fully algebraic treatment. Hence the standard Schwarzshild black hole can be obtained with $f^2(r)=1-{\gamma\over r}$ provided this is invertible. One can arrange this formally but  one can also proceed to address such issues using the topology of the field, for example by working over $k=\Q_p$ with $| \gamma/r|^p<1$ in the $p$-adic norm (i.e. a $p$-adic Schwarzschild black hole). One can also consider black holes over finite fields $\F_q$ but in this case the large kernel of $\extd$ if we use the usual differential calculus leads to solutions and phenomena that are artefacts of that; it would be interesting to consider instead reduced finite-dimensional versions of $\F_q[x_1,x_2,x_3]$ as our starting point and a more connected, noncommutative, differential calculus for this case, but following the pattern above. This will be considered elsewhere.

Finally, we compute the associated classical spacetime wave operator $\Delta$  in the setting of Theorem~4.7, for later reference. We  have now  also
\[ \(\bd t,\bd t\)=-{1\over f^2(r)},\quad \(\bd t,\bd r\)=\(\bd t,\omega_i\)=0\]
for the inverse spacetime metric. In our framework the spacetime wave operator on functions is $\bar\square=(\ ,\ )\bnabla\bd$ for the spacetime connection, exterior derivative and inverse metric.

\begin{corollary}\label{bsquare}The spacetime Laplace-Beltrami wave operator $\Delta$ associated to the metric in Theorem~4.7 is
\[\bar\square=-{1\over f^2}{\del^2\over\del t^2}+{1\over h^2}({2\over r}-{h'\over h}+{f'\over f }){\del\over\del r}+{1\over h^2}{\del^2\over\del r^2}+e_ie_i\]
(where we sum over $i$). \end{corollary}
\proof We first compute $\bar\square$,
\begin{eqnarray*} \bar\square \psi&=&\(\ ,\ )\bnabla\bd\psi=\(\ ,\ \)\bnabla(({\del\over\del t}\psi)\bd t+({\del\over\del r}\psi)\bd r+(e_i\psi)\omega_i)\nonumber\\
&=&(\ ,\ )(({\del^2\over\del t^2}\psi)\bd t\btens\bd t+{\del\over\del r}\psi)\bnabla\bd r+({\del^2\over\del r^2}\psi)\bd r\btens\bd r+(e_je_i\psi)\omega_j\btens\omega_i)\nonumber\\ 
&=&\left(-{1\over f^2}{\del^2\over\del t^2}+{1\over h^2}({2\over r}-{h'\over h}+{f'\over f }){\del\over\del r}+{1\over h^2}{\del^2\over\del r^2}+e_ie_i\right)\psi
\end{eqnarray*}
on a general function $\psi$ on $(\R^3\setminus \{0\})\times\R$. We showed only the terms in the outputs of  $\bnabla$ that are not immediately killed by the block-diagonal form of the inverse metric.  \endproof

Doing the same computation for the 3-geometry in Proposition~4.3 involves the same computations but without any of the terms involving $f$ and gives
\begin{equation}\label{bDelta} \bDelta_{LB} ={1\over h^2}({2\over r}-{h'\over h}){\del\over\del r}+{1\over h^2}{\del^2\over\del r^2}+e_ie_i\end{equation}
Note that the spatial part of $\bar\square$ differs from this by an extra ${f'\over fh^2}{\del\over\del r}$ and hence
\[ \bar\square = \beta{\del^2\over\del t^2}+\bDelta,\quad \bDelta=\bDelta_{LB}-{1\over 2}\beta^{-1}\bar g^{-1}(\bd\beta)\]
where $\beta=-1/f^2$, in accord with the general picture for this kind of metric explained in the proof of Corollary~3.4. Here  $\beta^{-1}\bd\beta=-{2 f'\over f}\bd r=-2\bar g({f'\over f  h^2}{\del\over\del r})$ using $\bar g$ as an operator. 

\subsection{Quantum black hole differential algebra} 

We are now ready to apply the machinery of Section 3 to spherically symmetric static spacetimes. Thus we let $(M,\bar g)$ be the classical Riemannian manifold in Proposition~4.3.  From Corollary~4.5 we have two types of conformal Killing vector and mainly we consider the first one and its divergence measure
\begin{equation}\label{tauh} \tau={r\over h(r)}{\del\over\del r},\quad \alpha={2\over h(r)}-1.\end{equation}
Thus the free function $h$ in the 3-geometry metric is now encoded in $\alpha$ as well as in the inverse 3-metric $\(\ ,\ \)$ and both are used in defining the calculus. We will not consider the noncommutative 4-geometry in any detail but from the classical wave operator in Corollary~4.8 compared with Corollary~3.4 we express the free function  $f(r)$  in the static metric in Proposition~4.6 and Theorem~4.7 as the free parameter
\begin{equation}\label{betaf} \beta=-{1\over c^2 f(r)^2}\end{equation}
of our noncommutative calculus. We have inserted the speed of light $c$ here. The functions $\mu,\nu$ there are now generically given by
\[ \mu(r)=-{1\over c^2 r^2}\int^r {h(r')\over f(r')^2}\, r'\extd r',\quad \nu(r)=e^{\int_1^r {h(r')-2\over r'}\extd r'}\int^r {e^{-\int_1^{r'} {h(r'')-2\over r''}\extd r''}\ h(r')\mu(r')\over r'}\extd r'\]
and provide the time part $2\Delta_0$ of the wave operator in the `finite difference' form in Proposition~3.3. It remains the case that Lemma~3.5 is a better route for its actual calculation. Then the wave operator according to Corollary~3.4 and Corollary~4.8 is
\begin{equation}\label{squareh}\square\psi(t)=2\Delta_0\psi(t)+(\left({1\over h^2}({2\over r}-{h'\over h}+{f'\over f}){\del\over\del r}+{1\over h^2}{\del^2\over\del r^2}+e_ie_i\right)\psi)(t+\lambda\alpha)\end{equation}
on normal ordered functions $\psi=\sum_n\psi_nt^n$.

Working through the constructions of Section~3 the quantum differential algebra associated to the classical geometry now takes the form:

\begin{proposition} The quantum calculus $\Omega^1(C(M)\rtimes\R)$ quantizing the classical picture in Proposition~4.6 with respect to  radial conformal Killing vector (\ref{tauh}) has relations
\[[x_i,x_j]=0,\quad [x_i,t]={\lambda\over h}x_i,\quad [\omega_i,x_j]=\lambda e_{ij}\theta',\quad [\extd r,x_i]={\lambda\over h^2}{x_i\over r}\theta',\quad [\theta',x_i]=0\]
\[ [\omega_i,t]=\lambda({1\over h}-1)\omega_i,\quad [\theta',t]=\lambda({2\over h}-1)\theta',\quad  [x_i,\extd t]=\lambda\extd x_i,\quad [\extd t,t]=\beta\lambda\theta'-\lambda\extd t.\]
\[ [\extd r,t]=\lambda(\extd({r\over h})-\extd r),\quad \extd g(r)=g'(r)\extd r+{\lambda\over 2 h^2}g''(r)\theta'\] 
for any function $g(r)$.\end{proposition}
\proof This is by application of Theorem~3.1 but can conveniently be found using the quantum form of the commutation relations (\ref{cr}). We use the same 1-forms $\omega_i$ as in Section~4.1  but note that these can also be viewed as the angular projection of the quantum differentials
\[   \omega_i=\bd x_i-{x_i\over r^2}x_j \bd x_j=\extd x_i-{x_i\over r^2}x_j \extd x_j=e_{ik}\extd x_k\]
so we need not  distinguish between the classical and quantum counterparts. To see this we use $\extd$ from Lemma~2.1 and  $\bDelta$ given below (\ref{bDelta}) to compute
\[ \extd x_i = \bd x_i + {\lambda\over 2 h^2}\left( {2\over r}(1-h^2))-{h'\over h}+{f'\over f})\right)\theta'.\]
One can then deduce properties of $\omega_i$ noting that $[e_{ij},t]=0$. We similarly have
\[ \extd r=\bd r+{\lambda\over 2 h^2}({2\over r}-{h'\over h}+{f'\over f})\theta'\]
and hence $\extd g=g'\bd r+{\lambda\over 2h^2}(({2\over r}-{h'\over h}+{f'\over f})g'+g'')\theta'$ comes out as stated. \endproof

One similarly has useful relations such as
\[ [\omega_i,r]=0,\quad[\extd r,g(r)]={\lambda g'(r)\over h^2}\theta',\quad [\extd r,x_i]=\lambda{x_i\over r h^2}\theta',\quad [\extd r,{x_i\over r}]=0,\quad x_i\omega_i=0\]
\[ \omega_i=\extd x_i-{x_i\over r}\extd r+{\lambda x_i\over h^2r^2}\theta',\quad [g(r),t]=\lambda {r\over h}g'(r),\quad [g(r),\extd t]=\lambda\extd g(r)\] 
on any function $g(r)$, where the last is mentioned for completeness (it applies to any function on $M$). Among our various relations we have  a closed algebra of $\extd r,\extd t,\theta'$ and functions of $r,t$. 

We now specialise to the case of a Schwarzschild black hole of mass $M$ where, from Theorem~4.7, we have classically $f=\sqrt{1-{\gamma\over r}}$ with $\gamma=2GM/c^2$ and $h=1/f$.  We regard $h,\tau,\alpha$ as associated to the 3-geometry and $\beta$ as a parameter in the calculus expressing the time component of the metric. We have 
\[ h={1\over \sqrt{1-{\gamma\over r}}},\quad \tau=r\sqrt{1-{\gamma\over r}}{\del\over\del r},\quad \alpha=2\sqrt{1-{\gamma\over r}}-1,\quad \beta=-{1\over c^2( 1-{\gamma\over r})}.\]
For these specific functions some of the formulae in Proposition~4.9 simplify, for example
\[ [\extd r,t]={\lambda\over 2}(\sqrt{h}-{1\over\sqrt{h}})^2\extd r-{\lambda^2\gamma^2\over 8 r^3}h\theta'.\]
In principle we can put the specific commutation relations into the theory in Section~3 to find the wave operator, obtained in practice say by Lemma~3.5. 

We also have a quantization by the Killing vector $\tau_3$ (say) in Corollary~4.5. Here
\begin{equation*} \tau_3={x_1 e_2-x_2 e_1\over h(r)^2},\quad \alpha=-1\end{equation*}
and we may take $h,\beta$ as desired, eg the Schwarzschild one. This of course breaks the rotational symmetry in the quantisation but may also be of interest as a complementary quantisation.

\section{Minimally coupled quantum black hole}

Although we have just quantised the Scwharzschild black hole in our approach, it remains the case that the resulting quantum differential algebra is extremely hard to compute with. In particular new methods will be needed to find the form of $\Delta_0$  in Proposition~3.5 because the function $\alpha$ enjoys rather poor commutation relations with $t$. In order to get a flavour of physical predictions we therefore consider now a simpler version in which we build the black hole by working in flat spacetime and afterwards allow for the spatial curvature `by hand' by  a process of minimal coupling. This has the merit of being fully computable and allows us our first quantum gravity predictions.

\subsection{Flat bicrossproduct model spacetime with gravity}

We apply the formalism of Sections~3 and~4 to $M=\R^3$ with flat metric with the conformal Killing vector given by  the degree operator $\rho$. Thus throughout this section,
\[ \tau=\rho=r{\del\over\del r},\quad \alpha=1\]
as a special case of Corollary~4.5 with $h(r)=1$. In this case  our calculus  $\Omega^1(C(\R^3)\rtimes\R)$ on generators becomes 
\[ [x_i,x_j]=0,\quad [x_i,t]=\lambda x_i,\quad [\extd x_i,x_j]=\lambda\delta_{ij}\theta',\quad [\theta',x_i]=0,\quad[\theta',t]=\lambda \theta'\]
\begin{equation}\label{bicrosscalc} [\extd x_i,t]=0,\quad [x_i,\extd t]=\lambda \extd x_i,\quad [\extd t,t]=\beta\lambda\theta'-\lambda\extd t.\end{equation}
Specialising the results of Section~4.2 or directly by similar methods,  we note

\begin{lemma} We  have $r\extd r=x_i\extd x_i+\lambda\theta'$ and a closed algebra of $\extd r,\theta',\extd t$ and functions $g$  of $r,t$ with
\[ \extd g(r)= g'(r)\extd r+ {\lambda\over 2}g''(r)\theta',\quad  [\extd r,g(r)]=\lambda g'(r)\theta',\quad [\theta',g(r)]=0,\quad   [\extd r,g(t)]=0\]
\[ [g(r),t]=\lambda r g'(r),\quad  [g(r),\extd t]=\lambda\extd g(r),\quad r g(t)=g(t+\lambda)r,\quad  \theta' g(t)=g(t+\lambda)\theta'\]
\end{lemma}
\proof All of these results follow from Section~4.2 setting $h=1$ or can be obtained directly by induction on the commutation relations on generators (some of which we have already seen in the proof of Proposition~3.3). One can also obtain them from Theorem~3.1 since this can be computed in any coordinates.  Note also that when $h=1$ we have $\extd x_i=\bd x_i$ and  $\extd r=\bd r+{\lambda\over r}\theta'$. \endproof 

Further relations useful in computations, all specialized from Section~4.2 or easily obtained directly, are
\[[\extd x_i, g(r)]=\lambda{x_i\over r}g'(r)\theta',\quad [\extd r,x_i]=\lambda{x_i\over r}\theta',\quad [\extd x_i,{x_j\over r}]=\lambda{e_{ij}\over r}\theta'\]
\[  x_i g(t)=g(t+\lambda)x_i,\quad  [\extd x_i,g(t)]=0,\quad [\extd r,{x_i\over r}]=0\]
\[ \omega_i=\extd x_i-{x_i\over r}\extd r+\lambda{x_i\over r^2}\theta',\quad [\omega_i,r]=0,\quad x_i\omega_i=0, \quad [\omega_i,x_j]=\lambda e_{ij}\theta',\quad [\omega_i,t]=0 \]
for any functions $g$. As we saw in the proof of Proposition~4.9 the $\omega_i$ coincide with their classical counterparts.

In our theory $\beta$ remains a free parameter of the differential calculus which by Corollary~3.4 one can interpret as the part of the quantum Laplacian that corresponds to the gravitational potential. Therefore the results of this section can be viewed as gravity bolted on to flat bicrossproduct spacetime. The physics of this will be justified elsewhere by taking the Newtonian gravity limit. The quantum spacetime wave operator from Corollary~3.4 has the form
 \begin{equation}\label{flatwavebeta} \square\psi(t)=2\Delta_0 \psi(t)+ \bDelta^{flat}\psi(t+\lambda)-{1\over 2\beta}\(\bd\beta,\bd\psi\)(t+\lambda)\end{equation}
on normal ordered spacetime functions $\psi=\sum_n \psi_n t^n$, where $\bDelta^{flat}$ is the classical flat space Laplacian on $\R^3$ and $\bDelta =\bDelta^{flat}-{1\over 2}\beta^{-1}\bd\beta$ is the spatial part of the classical Laplacian used in the construction.

We compute $\Delta_0$ and to keep things simple we limit ourselves to $\beta=\beta(r)$, i.e. to the spherically symmetric case. In this case the functions $\mu$ and $\nu$ in Proposition~3.3 are
\[ \mu={1\over r^2}\int^r  \beta(r') r'\bd r',\quad \nu={1\over r}\int^r  \mu(r') \bd r' .\]
where the constants of integration do not change the values of $\Delta_0$. This provides the `inner element' $\theta=\extd t-(\mu+\nu)\theta'$ and also in principle the time part of the geometry in the form of $\Delta_0$ according to Proposition~3.3.  From Lemma~3.5 we note that $\Delta_0$ depends linearly in $\beta$. Hence it is enough to compute it term by term for monomial $\beta$.

\begin{proposition}  For all real $m$ let $\beta={1\over r^m}$.
\begin{enumerate}
\item If $m=1$ 
\[ \mu={1\over r},\quad \nu= {\ln(r)\over r},\quad  \Delta_0g(t)={1\over r\lambda}({\del\over\del t}-\del_0)g(t+\lambda)\]
\item If $m=2$
\[ \mu={\ln(r)\over r^2},\quad \nu={1+\ln(r)\over r^2},\quad \Delta_0 g(t)={1\over r^2\lambda}\left(\del_0 g(t+2\lambda)-{\del\over\del t}g(t+\lambda)\right)\]
 \item If $m\ne 1,2$
 \[  \mu={1\over (2-m)r^m},\quad \nu={1\over (2-m)(1-m) r^m}\]
\[ \Delta_0 g(t)={1\over r^m}\left({g(t+\lambda)+(1-m)g(t-\lambda(1-m))-(2-m) g(t+\lambda m)\over \lambda^2(2-m)(1-m)}\right)\]
\item
\[ [\extd t, g(t)]+\lambda\del_0 g(t)\extd t={\lambda\over r^m}\left(\begin{cases} {\del\over\del t}g(t+\lambda)& m=2\\ {g(t+(m-1)\lambda)-g(t+\lambda)\over (m-2)\lambda} & m\ne 2\end{cases}\right)\theta'\]
\end{enumerate}
\end{proposition}
\proof The case $m=0$ is easily obtained directly or by noting that $\mu,\nu$ can be taken to be constants $\beta/2$ in Proposition~3.3. For the general case we integrate to find $\mu,\nu$ and use these in Proposition~3.3, or proceed along the same lines as next as for $m=1,2$. For $m=1$ we use one of the commutation relations in Lemma~5.1 in the form $g(t){1\over r}={1\over r}g(t+\lambda)$ to give us $\Delta_0 t^n$ as $1/r$ times
\[ \sum_{i=0}^{n-1} {(t+\lambda)^{n-1-i}-t^{n-1-i}\over\lambda}(t+\lambda)^i={n\over\lambda}(t+\lambda)^{n-1}-{(t+\lambda)^n-t^n\over\lambda^2}\]
which gives answer stated. For $m=2$ moving $1/r^2$ to the left gives us $1/r^2$ times
\[ \sum_{i=0}^{n-1} {(t+2\lambda)^{n-1-i}-(t+\lambda)^{n-1-i}\over\lambda}(t+\lambda)^i={(t+2\lambda)^n-(t+\lambda)^n\over\lambda^2}-{n\over\lambda}(t+\lambda)^{n-1}.\]
The generic case is similar if one wants to do it in the same way. Finally, the commutators $[\extd t,g(t)]$ are obtained from $\Delta_0$ (see Lemma~3.5) or from (\ref{dttn}). \endproof

The $m=0$ case $\beta=$constant has calculus
\[  [\extd t,g(t)]=-\lambda(\del_0g)\extd t+{\beta\over 2}(g(t+\lambda)-g(t-\lambda))\theta',\quad \extd g(t)=(\del_0g)\extd t+{\beta\lambda}(\Delta_0 g)\theta'\]
on functions $g(t)$, with $\Delta_0$ half the standard finite-difference Laplacian. In particular, we recover the standard $5$-dimensional version of the calculus cf \cite{Sit} on the bicrossproduct model or `$\kappa$-Minkowski' spacetime\cite{MaRue} in the precise conventions in which (in one space dimension lower) it appears  in 2+1 quantum gravity as a scaling limit of the standard quantum geometry of $C_q(SU_2)$ as this is stretched flat, see  \cite{MaSch,Ma:bsph}. We also recover
\[ \square_{\beta={\rm const}}\psi(t)={\psi(t+\lambda)+\psi(t-\lambda)-2\psi(t)\over\lambda^2}+ \bDelta^{flat}\psi(t+\lambda)\]
as familiar in the variable speed of light prediction \cite{AmeMa} for the bicrossproduct model flat spacetime without gravity. By contrast the choice $m=1$ results in the nonrelativistic limit in a Newtonian $1/r$ potential with noncommutative  deviations.

\subsection{Minimally coupled Schwarzschild black hole and predictions}

We are now ready to adapt the above to the black hole. First, we make the  particular choice
\begin{equation}\label{betaBH} \beta=-{1\over c^2(1-{\gamma\over r})}.\end{equation}
where $\gamma={2GM/c^2}$ will now be the Schwarzschild radius for a black hole of mass $M$. 
The Newtonian gravity point source mentioned above is the just the start of the geometric expansion of this $\beta$. We then have the differential calculus with the same relations as in Lemma~5.1 plus the remaining structure involving $\beta$  computed from Proposition~5.2.

Second, we `by hand' replace  $\bDelta^{flat}$ in (\ref{flatwavebeta}) by the Laplace-Beltrami operator $\bDelta_{LB}$ for the the specific 3-geometry in Proposition~4.3 that underlies the Schwarzschild black hole.  We also replace the flat inverse metric in the $\beta^{-1}\bd\beta$ term by the one for this 3-geometry. These steps are analogous to  working in flat 3-space coordinates followed by a process of `minimal coupling' where a covariant derivative is then put in by hand. Thus, we compute within the spatially flat space bicrossproduct model with the same $\Delta_0$ as above but adjust the wave operator from (\ref{flatwavebeta}) to
\[ \square_{BH}\psi(t)=2\Delta_0 \psi(t)+ \bDelta_{LB}\psi(t+\lambda)-{1\over 2\beta}\(\bd\beta,\bd\psi\)(t+\lambda)\]
on normal ordered spacetime functions $\psi=\sum_n \psi_n t^n$.  Explicitly,
\begin{equation}\label{BHwave}\square_{BH}\psi(t)=2\Delta_0\psi(t)+\left(({2\over r}-{\gamma \over r^2}){\del\over\del r}+(1-{\gamma\over r}){\del^2\over\del r^2}+e_ie_i\right)\psi(t+\lambda)\end{equation}
is our `minimally coupled' noncommutative black hole wave operator. 

It remains to study $\Delta_0$ further.  In order to effectively work with this we Fourier transform, i.e. consider the effect on functions with time dependence $\psi(t)=e^{\imath\omega t}$ where $\omega\in \R$ and we let $\lambda=\imath\lambda_p$.

\begin{proposition} For the Schwarzschild $\beta$ in (\ref{betaBH}) we have 
\[\Delta_0 e^{\imath \omega t}= {1\over c^2} D(\omega,r)e^{\imath\omega t}\]
where
\[ D(\omega,r)={1\over \lambda_p^2}\left({\sinh(\omega\lambda_p)}+e^{-\omega\lambda_p}(1-{\gamma\over r})\left(1-e^{\omega\lambda_p}-{\gamma\over r}\ln\left({e^{\omega\lambda_p}r-\gamma\over r-\gamma}\right)\right)\right)\]
has limits
\[ \lim_{\lambda_p\to 0}D(\omega,r)={\omega^2\over 2(1-{\gamma\over r})},\quad \lim_{r\to\infty}D(\omega,r)={\cosh(\omega\lambda_p)-1\over \lambda_p^2},\quad \lim_{r\to \gamma}D(\omega,r)={\sinh(\omega\lambda_p)\over \lambda_p^2}\]
\end{proposition}
\proof We do this by summing all the contributions in the geometric expansion of $\beta$ in the region $r>\gamma$ and using Proposition~5.2 for each term. Thus, setting $\zeta=e^{-\omega\lambda_p}$ for brevity,
\begin{eqnarray*}
-D(\omega,r)&=&{1\over 2\lambda^2}(\zeta+\zeta^{-1}-2)+{\zeta\gamma\over r\lambda}(\imath\omega-{(1-\zeta^{-1})\over\lambda})+{\zeta\gamma^2\over r^2\lambda}({\zeta-1\over\lambda}-\imath\omega)\\
&&+\sum_{m=3}^\infty{1\over r^m\lambda^2}\left({\zeta^m\over m-1}-{\zeta^{m-1}\over m-2}+{\zeta\over(m-1)(m-2)}\right)\\
&=&-{\zeta-\zeta^{-1}\over 2\lambda^2}+{1\over\lambda^2}(1-{\gamma\over r})\left(\zeta-1+{\imath\lambda\zeta\omega\gamma\over r}+{\zeta\gamma\over r}\ln\left({r-\gamma\over r-\zeta\gamma}\right)\right)
\end{eqnarray*}
which we write as stated. The limits are then easily obtained. For completeness, let us note that had we expanded the geometric series for $\beta$ appropriate to $r<\gamma$ we would have $\beta={1\over c^2}\sum_{m=1}^\infty ({r\over\gamma})^m$ and use Proposition~5.2 applied to $-m$, giving
\[ D(\omega,r)={1\over\lambda^2}\sum_{m=1}^\infty({r\over\gamma})^m\left({\zeta^{-(m+1)}\over m+2}-{\zeta^{-m}\over m+1}+{\zeta\over(m+1)(m+2)}\right)\]
which sums to the same expression as before. One can check that expanding the logarithm appropriately to $r$ small and $r$ large recovers the two different series. 
  \endproof

As the action of the finite difference $\del_0$ on $e^{\imath\omega t}$ is by $\del_0=(1-e^{\omega\lambda_p})/(\imath\lambda_p)$ and since all the frequency dependence of $D(\omega,r)$ is via $e^{\omega\lambda_p}$,  we can explicitly write the massless wave equation $\square\psi=0$ as
\begin{equation}\label{BHwave2}\left( {2\over c^2\lambda_p^2}{\mathcal D}(\imath\lambda_p\del_0)+({2\over r}-{\gamma \over r^2}){\del\over\del r}+(1-{\gamma\over r}){\del^2\over\del r^2}+e_ie_i\right)\psi=0\end{equation}
where ${\mathcal D}(X)= -X+{X^2\over 2}+(1-{\gamma\over r})(X-{\gamma\over r}\ln(1-{  X\over 1-{\gamma\over r}}))$.

Returning to Proposition~5.3, the first limit is correct as $2 D(\omega,r)/\omega^2=1/(1-{\gamma\over r})$ is then the classical coefficient in front of $-{1\over c^2}{\del^2\over\del t^2}$ in the wave operator. The second limit is also correct as $2 D(\omega,r)$ is then the time part of the Fourier transform of the wave operator in the flat bicrossproduct spacetime model as used in the VSL prediction\cite{AmeMa}. 
The third limit is more unexpected and we shall make some crude assumptions in order to get a `first impression' as to what this entails. 

{\em (a) Maximum  redshift and frequency dependence of the redshift}. For any static metric such as  $\beta^{-1}\bd t\btens\bd t+\bar g$ one has a standard argument for time dilation and associated frequency shift. Light emitted with frequency $\omega>0$ at, in our case, $r$ will appear at $\infty$ with frequency $\omega\sqrt{\beta(\infty)\over\beta(r)}$. By convention the redshift factor here is $1/(1+z)$. Thus for the classical black hole the redshift of an emission of frequency $\omega$ at $r$ will appear at $\infty$ with frequency $\omega/(1+z)$ where $(1+z)=1/\sqrt{1-{\gamma\over r}}$. We have not done an analysis of noncommutative photon propagation via $\square$ but as $2 D(\omega,r)/\omega^2$ enters in the same way as $\beta$ (namely in the $2\Delta_0$ part of $\square$)  we may expect that it has something like the same interpretation for each frequency mode. Thus 
\[ (1+z)=\sqrt{D(\omega,r)\over D(\omega,\infty)}\]
and for photons released closer and closer to the event horizon, the redshift factor $1+z$, which classically goes to infinity, we have a finite limit 
\[ (1+z)_{max}= \sqrt{\sinh(\omega\lambda_p)\over  \cosh(\omega\lambda_p)-1}.\]
 It should be stressed that this is not intended to give more than a  qualitative impression of the physics. 
For small $\omega\lambda_p$ we have $z_{max}\approx\sqrt{2\over\omega\lambda_p}$. For example, 
if $\lambda_p$ is Planck time and $\omega=10^{19}$ Hz (the upper end of the X-ray band) then  
\[ z_{max}\approx 5\times 10^{12}\]
but if the photon has planck scale energy-momentum then this maximum redshift tends to 1 as  $\omega\lambda_p\to \infty$. Thus the most energetic modes are not redshifted at all.

The frequency dependence of these formulae may  be a better route to detection. To assess this perhaps in the context of laser interferometry, consider a  laser source pointing away from the centre and consisting of a beam at frequency $\omega$ superimposed with a harmonic at some multiple of $n\omega$. However, on arrival at a distant receiver  they would no longer be the same multiple. Expanding 
\[ 2D(\omega,r)={\omega^2\over(1-{\gamma\over r})}\left(1-{2\over 3}{\omega\lambda_p\gamma\over r(1-{\gamma\over r})}+O((\omega\lambda_p)^2)\right)\]
we see that for small $\omega\lambda_p$, the harmonic will have smaller redshift factor than the base frequency and hence will appear to the distant observer as a little higher in frequency than the $n$'th harmonic.  Let $\omega'$ be the redshifted base frequency and $\omega''$ the redshifted harmonic. The deficit in distance per base cycle over which the harmonic completes its $n$ cycles is
\[{c\over\omega'}-n {c\over \omega''}={c\over\omega'}\left(1-\sqrt{D(n\omega,r)D(\omega,\infty)\over D(n\omega,\infty) D(\omega,r)}\right)\approx {(n-1)\over 3}{c\lambda_p\gamma\over r\sqrt{1-{\gamma\over r}}}\]
or approximately $n\gamma\over 3 r$ Planck lengths $\lambda_p$ error per base cycle on arrival. Taking a similar figure for the entire length $L$ of the journey (for our back-of-envelope estimate) we need
\[ L\sim {c^2\over \omega^2}{3  r \over n \gamma \lambda_p}\]
in order to accumulate one full cycle of phase error. For a 0.1 nanometer (X-ray) wavelength, $\gamma\over r$ around 0.1  (say), and $n=10$, we have some $L\sim 0.1$ light years which is modest by astronomical standards even if well beyond current reach. The figure would be worse using  infra red lasers but on the other hand it should not at all be necessary to accumulate a whole cycle of phase error to determine that $\omega''$ was not the same multiple of $\omega'$  and in that sense our preliminary estimate is very conservative.  Also note that we expect the frequency dependence of the redshift to apply to other gravitational potentials, not just to black holes, although clearly most of these would have an effective ${\gamma\over r}<<1$. 

{\em (b) Beckenstein-Hawking radiation.} This requires a certain amount of machinery to recompute from the noncommutative wave operator. However, at first sight the overall temperature to a distant observer should not change significantly for macroscopic  (non Planckian) black holes because the same factor in front of $-{1\over c^2}{\del^2\over\del t^2}$ enters into the computation of the acceleration and hence of the Unruh effect local  temperature near the horizon, which would also now be finite. The finiteness would appear to resolve the so-called `temperature paradox' whereby some authors have worried about the validity of the infinite temperature required at the horizon due to the infinite redshift from the horizon in the classical picture.  On the other hand, due to the frequency dependence of the redshift a black body spectrum at the horizon would no longer result in a black body after redshift. The more energetic modes should have less redshift thereby compressing the upper end of the distribution relative to the lower end.

 {\em (c) Wave operator at the horizon.} The limit $r\to \gamma$  in Proposition~5.3 and the limiting behaviour of the rest of the wave operator, means that the wave operator arbitrarily close to the event horizon in the standard Schwarzschild coordinates becomes
 \[ \lim_{r\to\gamma}\square_{BH}\psi(t)={\psi(t-\imath\lambda_p)-\psi(t+\imath\lambda_p)\over c^2 \lambda_p^2}+{1\over \gamma}{\del\over\del r}\psi(t+\imath\lambda_p)+e_ie_i\psi(t+\imath\lambda_p)\]
on normal ordered functions. We see that the singular $r-t$ sector of the classical wave operator drops down to what is conceptually a kind of `first order' differential operator as we approach $r=\gamma$ rather than blowing up in front of $\del^2\over\del t^2$ as it does classically. In a sense, the noncommutative deformation has smoothed out the classical coordinate singularity, at least as far as the wave operator is concerned.  Moreover, as the left hand side would be zero for a massless solution, one could think of this equation as a boundary condition for such solutions crossing the event horizon. Using the notation $\tilde\del_0\psi(t)={\psi(t)-\psi(t-2\lambda)\over2\lambda}$ for the finite difference (this is a version of $\del_0$ used elsewhere in the paper), and restricting for concreteness to $\psi$ a linear combination of the $Y^l_m$ spherical harmonics as regards angular dependence, we can write the condition  as
\[ {2\imath \over c}\tilde\del_0\psi= { \lambda_p\over  \gamma}{\del\over\del r}\psi- {\lambda_p\over\gamma^2}l(l+1)\psi\]
at  the horizon, where $\lambda_p$ is the Planck length. Assuming bounded spatial derivatives we see that in the classical limit where $\lambda_p\to 0$  or for infinitely large black holes as $\gamma\to \infty$, we will have  $\dot\psi=0$ at the horizon. However,  for a Planckian size black hole where $\gamma\sim \lambda_p$ we see that ${1\over c}\tilde\del_0\psi$ and ${\del\over\del r}\psi$ are comparable at least when $l=0$. 

Note also that for usual black holes the $r-t$ metric coefficients flip over in sign at the event horizon so that $r$ plays a role more like time inside the event horizon and vice-versa. In our case the function $D(\omega,r)$ while continuous in its real part at $r=\gamma$ acquires an imaginary part in a thin frequency-dependent layer at the horizon. If $\omega>0$ then this has thickness $\gamma (1-e^{-\omega \lambda_p})$ and is located just inside the classical event horizon,
\[\omega>0:\quad  \Im D(\omega,r)\ne 0,\quad \forall  r\in \gamma[e^{-\omega\lambda_p},1]\]
due to the negative argument of the logarithm.  One would need  artificially to use $\ln|\ |$ to avoid this `interregnum' layer just below the classical event horizon.  Below this layer, we have $D(\omega,r)$ negative as classically. Also note that as $r$ increases from below, the coordinate singularity is still present at the lower boundary $r=\gamma e^{-\omega\lambda_p}$ but is one degree lower so that $D(\omega,r)\sim\log(r-\gamma e^{-\omega\lambda_p})$ near this boundary, compared to classically.

When $\omega<0$ the picture is much the same except that the interregnum is reflected about the classical event horizon $r=\gamma$ and now lies just above it,
\[ \omega<0:\quad \Im D(\omega,r)\ne 0,\quad \forall  r\in \gamma[1, e^{-\omega\lambda_p}].\]
Thus the boundary  $r=\gamma e^{-\omega \lambda_p}$ of the interregnum where $D(\omega,r)$ is logarithmic now  lies just outside the the classical event horizon $r=\gamma$. Now $D(\omega,r)$ has a limit as $r\to\gamma$ when approached from below. We see that there appears to be an asymmetry in the treatment of positive and negative frequency modes. 

{\em (d) Singularity at the origin}. Finally, we note that $D(\omega,r)$ is again regular for small $r$ with expansion 
\[ D(\omega,r)=-{(\cosh(\omega\lambda_p)-1)(1+2e^{\omega\lambda_p})\over 3 \lambda_p^2\gamma }r + O(({r\over\gamma})^2)\]
deforming the classical behaviour but not too drastically for small $\omega\lambda_p$. However, for the Planckian velocities that might apply at the singularity at the origin, the effects appear to be similar to the well-known Planckian bounds at $r=\infty$. We recall that in the flat bicrossproduct spacetime model, the exponentially growing $\cosh(\omega\lambda_p)-1$ puts a bound on the spatial part of the wave operator. This does not imply but perhaps hints that some modes of the curvature might also be made finite, but this remains to be seen on a computation of more of the noncommutative geometry.

{\em (e) Numerical solutions}. Although an analytic study of the wave equation (\ref{BHwave2}) is beyond our scope here, we can get further  first impressions from a numerical study. As we have seen above, there are three regions of interest and we study them for given $\omega>0$ and Schwarzschild radius $\gamma=1$. One can think of the latter and $c=1$ as defining the reference scales in the system. For simplicity we limit ourselves to the sector with $l=0$ orbital angular momentum, i.e. without contribution from $e_ie_i$ in the wave operator. Computations have been done with MATHEMATICA.

\begin{figure}
\[(a)\includegraphics[scale=.9]{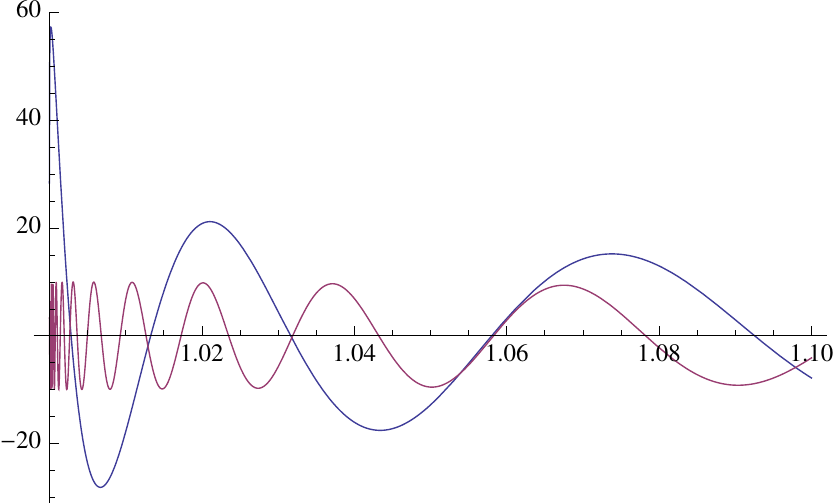}\]
\[(b) \includegraphics[scale=.9]{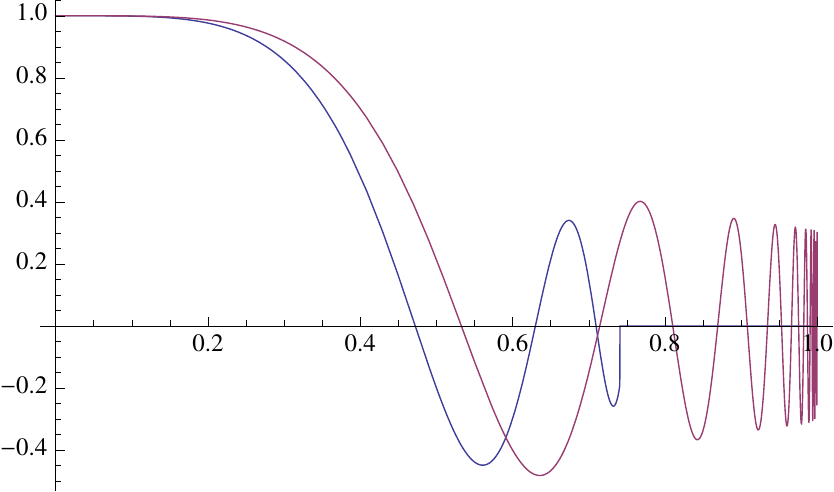}\]
\[ (c)\includegraphics[scale=.7]{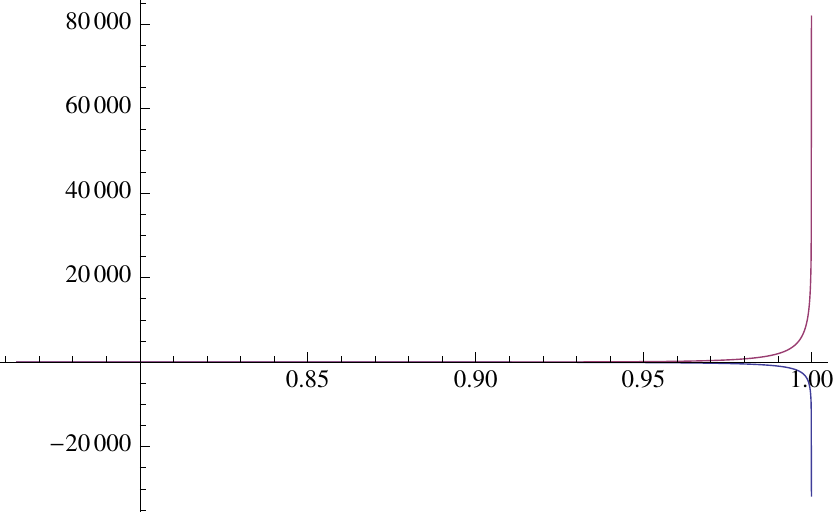}\includegraphics[scale=.7]{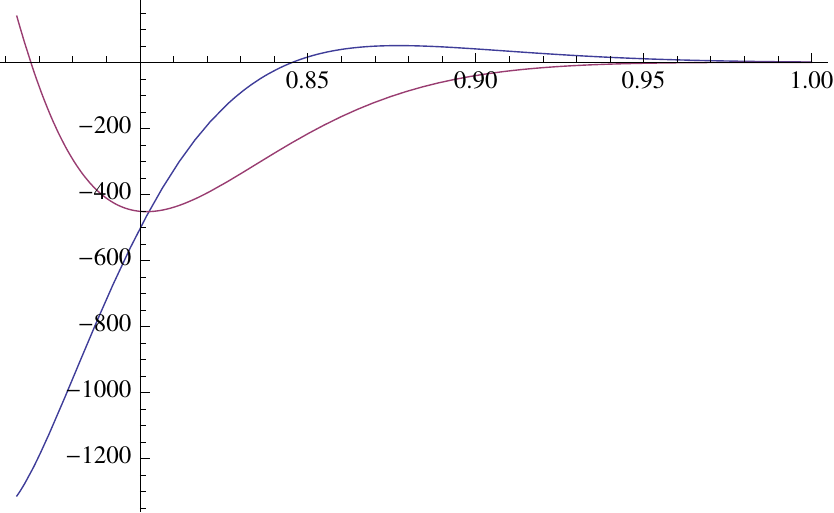}\]
\caption{Numerical solutions of the noncommutative wave equation on normal ordered functions with frequency $\omega>0$ and $\gamma=1$, and comparison with the classical black hole at same boundary conditions.  (a) Shows the exterior region $r>\gamma$ with waves appearing to have a finite frequency at the event horizon $r=\gamma$ as a new feature. (b) Shows the interior region $r<\gamma e^{-\omega\lambda_p}$ and the new possibility of standing waves with a finite number of `cycles'. The quantum solutions can be continued through from either side into (c) an interregnum region $\gamma e^{-\omega\lambda_p}\le r\le \gamma$ where they `amplify' and typically appear to diverge. The left plot shows solutions driven from the black hole interior and the right plot  from the black hole exterior. Shown are real and imaginary parts.}
\end{figure}

(i)  In this case the region $r>\gamma$  outside the event horizon is shown in Figure~1(a), comparing the classical case at $\lambda_p=0$ and the noncommutative case with the same boundary condition (namely, $\psi=1, \psi'=0$ far from the horizon at $r=10\gamma$). We used $\omega=10$ and $\lambda_p=1/10$ so that Planckian effects are pronounced. We see that whereas classical waves bunch up at the event horizon corresponding to an infinite redshift from the horizon to $r=\infty$, the noncommutative ones are not infinitely compressed at the event horizon. Although the solution shown appears regular at the event horizon one cannot be sure from the numerical solutions and a scaling limit analysis right at the event horizon suggests an unstable mode log-divergent at $r=\gamma$. One the other hand one can go the other way and set the amplitude and slope at or very close to the horizon as an alternative boundary condition that does not make sense classically. This supports the qualitative features discussed above concerning redshift and smoothing of the coordinate singularity.
 
(ii) The interior region $r<\gamma e^{-\omega \lambda_p}$ is shown in Figure~1(b), again comparing the classical and noncommutative cases with the same boundary conditions namely  $\psi=0$ and $\psi'$ fixed at $r=0$. In fact both classically and in the quantum case the solutions with such boundary conditions are highly insensitive to $\psi'$ at $r=0$ except in the immediate vicinity of $r=0$. We used $\omega=10$ and $\lambda_p=.03$ so that we start to see Planckian effects. Whereas classically the wave solutions bunch up infinitely at the boundary $r=\gamma$ of the black hole interior, the quantum case again appears to be regular at the boundary. This boundary in the quantum case is now located earlier (the interior region is strictly smaller than classically) but at that interior boundary one can achieve any chosen amplitude and slope, or indeed set these as boundary conditions at the interior boundary, something which makes no sense in the classical case. As well as the modes shown, there also appear to be wave modes with a log divergence at $r=0$ in both the classical and quantum cases more typically obtained by prescribing $\psi,\psi'$ at some point in the interior. 
 
In particular, we see from Figure~1(b) that as a purely quantum phenomenon we can now achieve `standing waves'  in which we can (if we want) have a certain number of `cycles', counting the descent from the initial start at $r=0$ as, say, 1/2 a cycle. Based on numerical work rather than analysis, we can adjust the frequency to get any (not necessarily whole) number of `cycles'  according mainly to the value of  $\omega\over\nu$, where $\nu={c\over\gamma}$ is the frequency associated to the Schwarzschild radius (something like the frequency of a wave spanning the black hole). For fixed $\omega\lambda_p$ the number of `cycles' appears to increase with $\omega\over\nu$ and for  $\omega\lambda_p=1$ one `cycle' is achieved at ${\omega\over \nu}=16$,  while for $\omega\lambda_p=10^{-5}$ one `cycle' is achieved at $\omega\lambda_p=0.45$.  This is a purely quantum phenomenon in the sense that to achieve a fixed pattern or number of cycles in the classical limit $\lambda_p\to 0$ we would appear to need an infinite $\omega$ and hence an infinitesimal size of black hole.  Moreover, for small $\omega\lambda_p$ the `cycles' become more and more skewed towards the boundary of the interior region. 
 
(iii) In between these two regions is the `interregnum' region $\gamma e^{-\omega \lambda_p}\le r\le \gamma$. The numerical solutions show this region to be unstable in the sense that for sufficiently high $\omega$ a prescribed non-zero amplitude or slope at either boundary appears to diverge or at least become very large when continued into the interregnum.  As for the interior region, the most relevant factor is  again $\omega\over\nu$, with divergence when this exceeds about 3 in the case of moderately  Plankian $\omega\lambda_p$ such as  $\omega=2.7,\lambda_p=1/10$ shown in Figure~1(c). In each case the boundary conditions are $\psi=1,\psi'=0$ at one edge (and real with real derivative).  For $\omega\lambda_p=1$ the corresponding critical values are around $3.5$ and 2.5 times $\omega\over\nu$, respectively, for boundary conditions set on the horizon side or the black hole interior side.   Other boundary conditions also appear to give solutions. The physical meaning of this interregnum region is unclear but may show a mechanism whereby an incoming external wave of sufficient $\omega$ is absorbed just below the event horizon rather than, as in the classical picture, being infinitely compressed just above it. Lower frequencies meanwhile appear to pass through the interregnum, possibly acquiring an imaginary component, and can be matched on the other side with wave solutions in the interior.

The case of $\omega<0$ is qualitatively much the same. This time the external noncommutative waves enter the interregnum before the classical event horizon at which classical solutions bunch up. Propagation inside the interregnum is  as in Figure~1(c) but with a left-right reflection of the figures.  The black hole interior is much the same as classically but with the noncommutative waves regular with the possibility of `standing waves'. Here the behaviour near the boundary of the black hole interior now at $r=\gamma$ is a left-right reflection of that for external waves approaching the event horizon in Figure~1(a). Thus  the positive and negative frequencies behave similarly but a little differently. Clearly, all of these matters require further investigation.


\begin{thebibliography}{99}

\bibitem{AmeMa} G. Amelino-Camelia \& S. Majid, Waves on noncommutative spacetime and gamma-ray bursts, Int. J. Mod. Phys. A15 (2000) 4301-4323

\bibitem{BegMa1}
E.J. Beggs \& S. Majid, Semiclassical differential structures, Pac. J. Math.224 (2006) 1-44
\bibitem{BegMa2}
E.J. Beggs \& S. Majid, Quantization by cochain twists and nonassociative differentials,  J. Math. Phys., 51 (2010) 053522 (32pp)
\bibitem{BegMa4}
E.J. Beggs \& S. Majid, *-Compatible connections in noncommutative Riemannian geometry,   J. Geom. Phys. (2011) 95--124
\bibitem{BegMa3}
 E.J. Beggs \& S. Majid, Nonassociative Riemannian geometry by twisting, J. Phys. Conf. Ser. 254 (2010) 012002 (29pp)


\bibitem{Con}
A. Connes, Noncommutative Geometry, Academic Press (1994).

\bibitem{MH}
A. Dimakis \& F. M\"uller-Hoissen, A noncommutative differential calculus and its relation to gauge theory and gravitation, Int. J. Mod. Phys. A (Proc. Suppl.) Vol. no. 3A, (1993) 474--477; Noncommutative differential calculus, gauge theory and gravitation, Gottingen preprint 1992

\bibitem{DV1}M. Dubois-Violette \& T. Masson, On the first-order operators in bimodules, Lett. Math. Phys. 37 (1996) 467--474.
\bibitem{DV2}M. Dubois-Violette \& P.W. Michor, Connections on central bimodules in noncommutative differential geometry, J. Geom. Phys. 20 (1996) 218 --232

\bibitem{FreMa}
L. Freidel and S. Majid, Noncommutative harmonic analysis, sampling theory and the Duflo map in 2+1 quantum gravity, Class. Quant. Gravity 25 (2008) 045006 (37pp)



\bibitem{Ma:pla}
S. Majid, Hopf algebras for physics at the Planck scale, J. Class. Quant. Gravity 5 (1988) 1587-1607


\bibitem{MaRue} S. Majid \& H. Ruegg,  Bicrossproduct structure of the $\kappa$-Poincare group and non-commutative geometry, Phys. Lett. B. 334 (1994) 348-354


\bibitem{Ma:rie}S. Majid, Quantum and braided group Riemannian geometry, J. Geom. Phys. 30 (1999) 113-146

\bibitem{Ma:rieq}S. Majid, Riemannian geometry of quantum groups and finite groups with nonuniversal differentials, Commun. Math. Phys. 225 (2002) 131-170


\bibitem{Ma:spo} S. Majid, Noncommutative model with spontaneous time generation and Planckian bound, J. Math. Phys. 46 (2005) 103520 (18pp)


\bibitem{Ma:bsph} S. Majid, $q$-Fuzzy spheres and quantum differentials on $B_q[SU_2]$ and $U_q(su_2)$, Lett. Math. Phys. (2011) 21pp

\bibitem{MaSch}
S. Majid \& B. Schroers, q-Deformation and semidualisation in 3D quantum gravity, J. Phys
A 42 (2009) 425402 (40pp)


\bibitem{Mou}J. Mourad, Linear connections in noncommutative geometry, Class. Quantum Grav. 12 (1995) 965 -- 974

\bibitem{Sit} A. Sitarz, Noncommutative differential calculus on the $\kappa$-Minkowski space, Phys.Lett. B 349
(1995) 42-48


\bibitem{Schupp} P. Schupp \& S.N. Solodukhin,  Exact black hole solutions in noncommutative gravity, arXiv:0906.2724[hep-th]


\end{thebibliography}
\end{document}